\documentclass[a4paper,12pt]{article}
\usepackage{indentfirst}
\usepackage{amsmath}
\usepackage{amsthm}
\usepackage{amssymb}
\usepackage{tikz}
\usepackage{fullpage}
\usepackage{mathrsfs}
\usepackage{enumitem}
\usepackage{hyperref}
\usetikzlibrary{decorations.markings}
\usetikzlibrary{cd}

\title{Bounding Projective Hypersurface Singularities}
\author{Ben Castor}
\date{}

\theoremstyle{plain}
\newtheorem{thm}{Theorem}[section]
\newtheorem{cor}[thm]{Corollary}

\newtheorem{lemma}[thm]{Lemma}
\newtheorem{prop} [thm]{Proposition}

\theoremstyle{definition}
\newtheorem{exmp}[thm]{Example}

\begin{document}
\maketitle

\begin{abstract}
We compare several different methods involving Hodge-theoretic spectra of singularities which produce constraints on the number and type of isolated singularities on a projective hypersurface of fixed degree.  In particular, we introduce a method based on the spectrum of the nonisolated singularity at the origin of the affine cone on such a hypersurface, and relate the resulting explicit formula to Varchenko's bound. 
\end{abstract}

\section*{Introduction}

The spectrum developed by Steenbrink \cite{St85} has served as an invaluable tool in understanding the monodromy about complex singularities, while providing a powerful and easily computable invariant for isolated singularities.  It is obtained from a natural mixed Hodge structure on the cohomology of the Milnor fiber of such a singularity, by combining information on the Hodge filtration and eigenvalues of the semisimple part of monodromy. The Saito-Steenbrink formula (conjectured in \cite{St89}, proved in \cite{Sa91}) expands this theory by relating the spectrum of a \emph{nonisolated} singularity to that of the isolated singularity obtained via its Yomdin deformation by a power of a linear form. 

We arrive in this setting by taking the affine cone on a projective hypersurface with isolated singularities, and a separate formula from \cite{St89} provides the contributions of the isolated singularities to the spectrum of the (nonisolated) cone point.  In the present paper, we show how to combine these results to bound the number of singularities of a particular type (or combination of types) that can be present on the hypersurface.  We also reduce the computations of the method down to a bound similar to that of Varchenko \cite{Var83}, and conjecture that this bound is the same.  

A projective hypersurface $X \subset \mathbb{P}^n$ of degree $d$ will thus be the main object under discussion.  In Section \ref{S1}, we derive a purely combinatorial formula for the Hodge numbers of $X$.  Then in Sections \ref{S2}-\ref{S5} we review several attempts to use Hodge-theoretic methods to bound the number of singularities of any particular deformation class on $X$ in terms of only $n$ and $d$. Section \ref{S2} concerns itself with the vanishing cycle sequence method, which utilizes only the exactness of the vanishing cycle sequence, and the properties of cohomological objects, summarized in Prop. \ref{prop2.1}.  However, this method has some limitations, and produces no bound for $A_1$ singularities when $n$ is even (i.e. when $X$ has odd dimension).

In Section \ref{S3}, we review Varchenko's bounding method \cite{Var83}, which was recently revisited in a beautiful expository article by van Straten \cite{vS20}. In Section \ref{S4}, we restate and give an initial generalization of the ``conical method'' worked out by Steenbrink and T. de Jong (cf. \cite{St89}) for bounding the number of double points on $X$, which duplicated Varchenko's bound in that case.   The bound is obtained by combining a separate formula of Steenbrink (given here as Theorem \ref{th4.2}) with his conjecture in the case (known to him) where the isolated singularities on $X$ are of Pham-Brieskorn type.  Steenbrink conjectures in \cite{St89} that a version of this latter formula can be proven in much higher generality, which Saito did in \cite{Sa91}. 

Siersma \cite{Si90} proved a weaker version of Steenbrink's conjecture by focusing solely on the characteristic polynomials of the monodromy operators (thus dropping the information on the Hodge filtration contained in the spectrum).  An attempt to bound the number of singularities based on Siersma's work is given in Section \ref{S5}.  Though sometimes better than the ``naive'' bound from Section \ref{S2}, this bound is still relatively weak; on the other hand, Siersma's results allow us to compute the vertical monodromies required for the more general Saito-Steenbrink formula to work.  

In section \ref{S6}, we use this information to generalize the main results of Section \ref{S4} to arbitrary isolated hypersurface singularities, and we show that in this case the power of the general linear form in the Yomdin deformation of the cone can be chosen to be any $k>d$.  This results in our Theorem \ref{thm1}, which we then make completely explicit in Theorems \ref{th7.1} and \ref{th7.7}.  We conclude the paper by showing that Varchenko's bound, while less discrete, implies Theorem \ref{th7.7}, but we may add another set of discrete bounds, which, in practice, are often extraneous, to arrive at a discrete bound stronger than Varchenko's.
 
The main result of this paper is given in Theorem \ref{thm1} which proceeds as follows; let $f$ be a homogeneous polynomial of degree $d$, with $X:=\widetilde{V}(f) \subseteq \mathbb{P}^n$ having only finitely many isolated singularities $P_i$. We write $\sigma_{g_i} = \sum_{j=1}^{\mu_i} [\lambda_{i,j}]$ for their spectra, and set $\alpha_{i,j} = d\lambda_{i.j} - \lfloor d \lambda_{i,j} \rfloor$.  Then for a suitably general linear form $\ell$, and any $k>d$, $f_k =f + \epsilon \ell^k$ has an isolated singularity at zero, and
\begin{equation}\label{eq1} 
\sigma_{f_{k},0}= \gamma_d^{*(n+1)} - \sum_{i,j} \bigg[\lambda_{ij}-\frac{\alpha_{ij}}{d} \bigg]* \beta_{d} + \sum_{i,j} \bigg[\lambda_{ij}-\frac{\alpha_{ij}}{k} \bigg]* \beta_{k},
\end{equation}
where $\gamma_d^{*(n+1)}=(\sum_{m=1}^{d-1}[-\tfrac{m}{d}])^{*(n+1)}$ is the spectrum of $\sum_{k=0}^n x_k^d$.

In particular, the effectiveness of \eqref{eq1} implies that of the right as an element of the free abelian group $\mathbb{Z}[\mathbb{Q}]$.  This leads to our ``conical bound'', which is given in a reduced form in Theorem \ref{th7.7}, for isolated singularities on $X$.  Here are a few concrete examples:
\begin{itemize}[leftmargin=0.5cm]
\item For $X \subset \mathbb{P}^3$ of degree $d$ with only $n_6$, $\widetilde{E_6}$, $n_7$ $\widetilde{E_7}$, and  $n_8$ $\widetilde{E_8}$ singularities, the sum $7n_6 +8n_7+9n_8$ is bounded by the polynomial $$\frac{277}{432}d^3 - \frac{23}{36}d^2 +\frac{53}{12}d + \frac{1}{2},$$ and in the case of only $\widetilde{E_6}$ singularities we have $$n_6 \leq \frac{31}{378}d^3 - \frac{13}{126}d^2 +\frac{4}{7}d + \frac{1}{14}.$$ This compares favorably to the ``naive'' Hodge theoretic bound $$n_6 \leq \frac{1}{9}d^3-\frac{1}{3}d^2+\frac{7}{18}d- \frac{1}{6}$$ resulting from the vanishing cycle sequence in Section \ref{S2}.
\item If $X \subseteq \mathbb{P}^4$ has degree $d$ with only $A_{2m+1}$ singularities, then the number of these is bounded by
     $$ r \leq  \begin{cases} 
     \frac{1}{2m+1}\left[ \frac{115}{192}d^4 -\frac{115}{48}d^3+ \frac{185}{48}d^2 - \frac{35}{12}d + 1 \right]  & d\equiv 0 \mod 2 \\
        \frac{1}{2m+1}\left[ \frac{115}{192}d^4 -\frac{115}{48}d^3 + \frac{355}{96}d^2 - \frac{125}{48}d + \frac{45}{64} \right]  & d\equiv 1 \mod 2, d>m+1\\
   \end{cases}
    $$
This easily beats the naive bound, which is asymptotic to $\tfrac{11}{24m}d^4$.
\item Many of the surface singularities ``with $K3$ tail'' classified in \cite[\S3]{KL2}, which include for instance the Dolgachev singularities, can occur on a quartic hypersurface in $\mathbb{P}^3$  (Singularities whose simplest form involves powers greater than $4$ can have analytically equivalent forms where this is not the case.)  The 3rd entry in [op. cit., Table 2], given by $x^2+y^6+z^6$, is not ruled out by the ``naive'' bound; but it is prohibited by the conical bound.
\end{itemize}
An appendix (Section \ref{S8}) contains tables providing additional comparisons of the various bounds.
  
\subsection*{Acknowledgements}
I would 
like to thank my Ph.D. advisor Matt Kerr, who advised me on the material for tackling this paper and was generously receptive to my persistent questioning on the material.

\section{A Formula for Hodge Numbers of Smooth Projective Hypersurfaces}\label{S1}

In this section, we derive explicit formulas for the Hodge numbers $h^{i,j}(H^{n-1}(X))$ of a smooth projective hypersurface $X \subseteq \mathbb{P}^n$ of degree $d$ in terms of $n$ and $d$.  This is of course very classical and included mainly for reference.  Since these numbers are indepedent of $X$, we denote them by $h^{i,j}_{n,d}$.  We will also write 
$$  [h_{n,d}^{i,j}]' =  h_{n,d}^{i,j} - \delta_{i,j} $$
for the primitive Hodge numbers (denoted $h_0^{i,j}(X)$ in \cite{Ar}).  It is well known that $h^{i,j}_{n,d} = \delta_{i,j}$ (Kronecker delta) for $ n-1 \neq i + j \leq 2n-2$. However, the middle Hodge numbers $h^{k,n-1-k}_{n,d}$ are much more complicated to calculate. 

\begin{thm}\label{th1.1}
For a smooth hypersurface $X \subseteq \mathbb{P}^{n}$ of degree $d$, the middle primitive Hodge numbers $[h^{k,n-1-k}_{n,d}]'$ (where $k \leq \frac{n-1}{2}$) are given by:
 $$[h^{k,n-1-k}_{n,d}]' =  (-1)^{n}\sum_{i=0}^k  (-1)^i\binom{n+1}{i}\binom{n-(k+1)d+(d-1)i}{n} $$
In particular, if $d>n$, then:
 $$ [h^{k,n-1-k}_{n,d}]' = \sum_{i=0}^{k} (-1)^{i}\binom{n+1}{i} \binom{(k+1)d-1-(d-1)i}{n}. $$
    
\end{thm}

The proof makes use of the following

\begin{lemma}
We let $\chi$ denote the Euler characteristic in the context defined in \cite{Ar}. Let $X \subseteq \mathbb{P}^n$ be a  smooth hypersurface of degree $d$. Recall, $\Omega_X^p = \Omega^p_X(0)$ and $\Omega^0_X(i)= \mathcal{O}_X(i)$ by definition. Then:

    \begin{enumerate}
        \item $\chi(\Omega_X^0(i))=\binom{i+n}{n}-\binom{i+n-d}{n}$
        \item $\chi(\Omega_{\mathbb{P}}^k(i))=\sum_{j=0}^k(-1)^j \binom{n+1}{k-j}\binom{n-k+i+j}{n}$
        \item $P^k_d(i)= \chi(\Omega^k_{\mathbb{P}}(i)) - \chi(\Omega^k_{\mathbb{P}}(i-d)) = \sum_{j=0}^k (-1)^j \binom{n+1}{k-j}\big[\binom{n-k+i+j}{n}-\binom{n-k-d+i+j}{n} \big]  $
        \item $\chi(\Omega^k_X(i))=P^k_d(i)-\chi(\Omega^{k-1}_X(i-d))$
        \item $\chi(\Omega^k_X)= \sum_{m=0}^k(-1)^{k-m}\sum_{j=0}^m (-1)^j\binom{n+1}{m-j}\big[\binom{n-m-(k-m)d+j}{n}-\binom{n-m-(k-m+1)d+j}{n} \big] $
    \end{enumerate}

\end{lemma}

\begin{proof}
     (1)-(4) of the lemma are just \cite[Prop 17.3.2]{Ar} stated in such a way that it is easier to follow the steps of the recurrence relation. 
     
     For the proof of (5), note that by (1), $$\chi\left(\Omega_X^{k-k}(-kd)\right)= \chi\left(\Omega_X^0(-kd)\right)= \binom{-kd+n}{n}-\binom{-kd+n-d}{n}$$
     And by (4):
     $$ \chi \left(\Omega_X^{k-(k-m)}(-(k-m)d \right) = P^{k-(k-m)}_d\left(-(k-m)d \right) - \chi\left(\Omega_X^{k-(k-m+1)}(-(k-m+1)d)\right)  $$
    These two facts together inductively give us:
    $$ \chi\left( \Omega^k_X \right) = \Bigg[\sum_{m=1}^k (-1)^{k-m} P^m_d(-(k-m)d) \Bigg]+(-1)^k\Bigg[ \binom{n-kd}{n} - \binom{n-(k+1)d}{n} \Bigg]$$
    And substituting (3) for $P^m_d(-(k-m)d)$, and noting that the second summand is just the case $m=0$, we get (5). 
\end{proof}

Now we use the last part of our lemma and the relationship between $[h^{k,n-1-k}_{n,d}]'$ and $\chi(\Omega^k_X)$ to prove our formula.

\begin{proof} (of Theorem)
     We know from \cite[Lemma 17.3.1]{Ar} that 
     $$[h^{k,n-1-k}_{n,d}]'= (-1)^{n-1-k}\chi(\Omega^k_X)+(-1)^n$$
     Using (5) in the above lemma we easily get a sum for $(-1)^{n-1-k}\chi(\Omega^k_X)$, setting $m-j=i$, and noting that $\binom{n-i}{n}=0$ for $1\leq i \leq k$, we get: 
    $$[h^{k,n-1-k}_{n,d}]' =  (-1)^{n}\sum_{i=0}^k \bigg[ (-1)^i\binom{n+1}{i}\binom{n-(k+1)d+(d-1)i}{n} \bigg]$$
    We know that for $q>0, \quad p<0$ $\binom{p}{q}=(-1)^q\binom{q-p-1}{q}$. 
    If $d>n$, Then  \linebreak
    $n-(k+1)d+(d-1)i = (n-d)- kd + (d-1)i  \leq -1 -kd + (d-1)k \leq -1 < 0 $ for $0\leq i \leq k$.
    Therefore:
    $$ [h^{k,n-1-k}_{n,d}]' = \sum_{i=0}^{k} (-1)^{i}\binom{n+1}{i} \binom{(k+1)d-1-(d-1)i}{n}. $$
\end{proof}

\section{The Vanishing Cycle Sequence Method}\label{S2}

Let $\pi\colon \mathcal{X}\to \Delta$ be a family of projective hypersurfaces $X_t:=\pi^{-1}(t)\subset \mathbb{P}^n$ over a disk about $t=0$, with smooth total space.  We assume that the fibers over $\Delta^*:=\Delta\setminus\{0\}$ are smooth of degree $d$, and $X_0$ has only isolated singularities.  We write $f\in \mathbb{C}[x_0,\ldots,x_n]$ for the homogeneous polynomial of degree $d$ cutting out $X_0\,(=\widetilde{V}(f))$.  In this setting, the vanishing cycle sequence is an exact sequence of mixed Hodge structures of the form 

 \begin{center}
    \begin{tikzcd}[column sep=scriptsize]
        0  \rar &  H^{n-1}(X_0) \rar & H_{lim}^{n-1}(X_t) \rar & H_{van}^{n-1}(X_t) \rar["\delta"] & H^{n}(X_0) \rar & H_{lim}^{n}(X_t) \rar & 0.
    \end{tikzcd}
\end{center}
(The mixed Hodge structures are induced by the nearby cycles triangle in the derived category of mixed Hodge modules on $X_0$, cf. \cite{KL1}.)  In particular, when $n$ is odd, $H_{lim}^{n}(X_t)=0$ and the sequence simplifies to

\begin{center}
    \begin{tikzcd}[column sep=scriptsize]
        0  \rar &  H^{n-1}(X_0) \rar & H_{lim}^{n-1}(X_t) \rar & H_{van}^{n-1}(X_t) \rar["\delta"] & H^{n}(X_0) \rar & 0
    \end{tikzcd}
\end{center}
We will use these sequences and the following facts to induce an inequality bounding the number and type of singularities of $f$.

\begin{prop} \label{prop2.1}
    \begin{enumerate}[label=(\alph*)]
        \item $H^{n-1}(X_t)\neq 0$ has Hodge numbers $h^{p,n-p-1}(H^{n-1}(X_t))=h^{p,n-p-1}_{n,d}$.
        \item Suppose $X_0$ has isolated singularities $p_1, \ldots p_r$ given locally by polynomials $g_1,\ldots,g_r$ (in $n$ variables), with Milnor fibers $Y_{g_i}$. Then we have an isomorphism $H^{n-1}_{van}(X_t)\cong \bigoplus_i \tilde{H}^{n-1}(Y_{g_i})$, where the MHSs on the reduced Milnor fiber cohomologies again come from Saito's theory.
        \item $\delta$ is a map of pure weight $n$ and, as such, can only be nonzero on $(p,q)$ parts for $p+q=n$.
        \item The Hodge numbers $h_{lim}^{p,q}$ of $H_{lim}^{n-1}(X_t)$ satisfy $\sum_q h_{lim}^{p,q} = h^{p,n-p-1}_{n,d}$ for each fixed $p$.  Moreover, they are symmetric about the lines $p=q$ and $p+q=n-1$.
        \item Suppose $U$ is a complex algebraic variety of dimension $m$. Then the values of $(p,q)$ for which the Hodge numbers $h^{p,q}(H^k(U))\neq 0$ satisfy:
            \begin{enumerate}
            \item $0 \leq p,q \leq k$;
            \item if $k>m$ then $k-m \leq p,q \leq m$;
            \item if $U$ is smooth then $p+q \geq k$; and
            \item if $U$ is compact then $p+q \leq k$.
            \end{enumerate}        
        \item With $X_0$ the singular fiber above, if $n={2m+1}$ is odd,  we have $h^{m,m}(H^{2m}(X_0)) \geq 1$
    \end{enumerate}
\end{prop}

\begin{proof}
    \begin{enumerate}[label=(\alph*)]
        \item This follows from the fact that $X_t \subset \mathbb{P}^{n}$ is a smooth hypersurface of degree d and complex dimension $n-1$, and therefore has a pure Hodge decomposition. 
        \item This follows from \cite[Theorem 5.44]{SP}.
        \item See \cite[Prop. 5.5]{KL1}.
        \item By \cite[p263,285]{SP} $\dim F^mH^k(X_t)= \dim F^m H^k_{lim}(X_t) $ for any $k,m \in \mathbb{N}$, where Steenbrink denotes $H^k_{lim}(X_t)= H^k(X_{\infty})$.  So $\sum_q h^{p,q}_{lim}=\dim (\mathrm{Gr}_F^p H^{n-1}_{lim}(X_t))$ must equal $\dim(\mathrm{Gr}_F^p H^{n-1}(X_t))=h^{p,n-p-1}_{n,d}$.
        \item This is just \cite[Theorem 5.39]{SP}
        \item Let $X_0 \subset \mathbb{P}^{2m+1}$, and let $\widetilde{X_0}$ denote a resolution of singularities:  pictorially,
            \begin{center}
            \begin{tikzcd}
                 X_0 \arrow[r, hook, "\imath"]  & \mathbb{P}^{2m+1}\\
                 \widetilde{X_0} \arrow[u, two heads, "\pi"] \arrow[ur, "\tilde{\imath}"']
            \end{tikzcd}
            \end{center}
        where $\imath$ and $\pi$ are the usual inclusion and projection maps.  This produces a commutative diagram of MHSs
        \begin{center}
            \begin{tikzcd}
                H^k(\mathbb{P}^{2m+1}) \arrow[r,"\imath^{*}"] & H^k(X_0) \arrow[rr,"\imath_{*}"] \arrow[dr,"\pi^{*}"'] & & H^{k+2}(\mathbb{P}^{2m+1})(1) \\
                && H^k(\widetilde{X_0}) \arrow[ur,"\tilde{\imath}_{*}"']
            \end{tikzcd}
        \end{center}    
where $\tilde{\imath}_{*}$ is the Gysin map and $\imath_*$ was defined by the composition $\tilde{\imath}_*\circ\pi^*$.  However, the map $\imath_*\circ \imath^*$ is really just the Lefschetz operator $L^k: H^k(\mathbb{P}^{2m+1})\rightarrow H^{k+2}(\mathbb{P}^{2m+1})$ which is an isomorphism (given by cupping with a hyperplane class $[H]$) for all $0\leq k\leq 4m$ by Hard Lefschetz. This implies that $\imath_*$ is surjective.  Taking $k=2m$, its image has type $(m,m)$, and so $(H^2m(X_0))^{m,m}$ cannot be zero.\end{enumerate}
\end{proof}

\begin{prop}\label{prop2.2}

    Let $f \in \mathbb{C}[x_0, x_1, x_2, x_3]$ be homogeneous polynomial of degree $d > 3$, and let $\widetilde{V}(f) \subseteq \mathbb{P}^3$ have only isolated $\widetilde{E_6}$ singularities. Then the number $r$ of singular points is bounded by
    $$ r \leq \frac{1}{6} (h^{1,1}_{3,d} - 1) = \frac{1}{6} \left[ \binom{2d-1}{3} - 4 \binom{d}{3} \right] = \frac{1}{9}d^3-\frac{1}{3}d^2+\frac{7}{18}d- \frac{1}{6}. $$

\end{prop}

\begin{proof}
    Since n=3 is odd, we have the following exact sequence of mixed Hodge structures:
    
    \begin{center}
        \begin{tikzcd}[column sep=scriptsize]
            0  \rar &  H^{2}(X_0) \rar & H_{lim}^{2}(X_t) \rar & H_{van}^2(X_t) \rar["\delta"] & H^{3}(X_0) \rar & 0
        \end{tikzcd}
    \end{center}
    We may visualize the exact sequence in terms of the mixed Hodge numbers using the Hodge-Deligne diagrams, where $W_k$ denotes the vector subspaces with weight $k$:

    
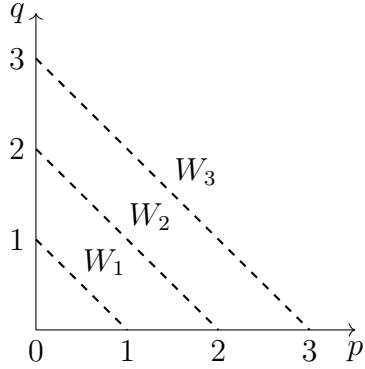
\vspace{2mm}
    
    \begin{center}
    \begin{tikzpicture}[scale=.6]

        \draw[->] (0,0) -- (7,0) node[anchor=north] {$p$};
        \draw	(0,0) node[anchor=north] {0}
        		(2,0) node[anchor=north] {1}
        		(4,0) node[anchor=north] {2}
        		(6,0) node[anchor=north] {3};
        
        \draw[->] (0,0) -- (0,7) node[anchor=east] {$q$};
        \draw	(0,2) node[anchor=east] {1}
        		(0,4) node[anchor=east] {2}
        		(0,6) node[anchor=east] {3};
        
        \draw[thick,dashed] (0,2) -- (2,0)
                            (0,4) -- (4,0)
                            (0,6) -- (6,0);
        \draw   (1.5,1.5) node {$W_1$} 
                (2.5,2.5) node {$W_2$} 
                (3.5,3.5) node {$W_3$}; 

    \end{tikzpicture}
    \end{center}

    Since $X_0$ is a compact variety of dimension $2$, we know that $H^2(X_0)$ has the form below by Prop. \ref{prop2.1}[(e),(f)]. Additionally, $H^2_{van}(X_t)$ has Hodge numbers given by $r$ times the Hodge numbers of $H^2(Y_g)$, where $g: \mathbb{C}^3 \rightarrow \mathbb{C}$ has a single $\widetilde{E_6}$ singularity. These are calculated explicitly in \cite{St77} as $h^{1,2}(H^2(Y_g))=h^{2,1}(H^2(Y_g))=1$ and $h^{1,1}(H^2(Y_g))=6$, and $0$ for the rest, giving the diagram below:
    
    
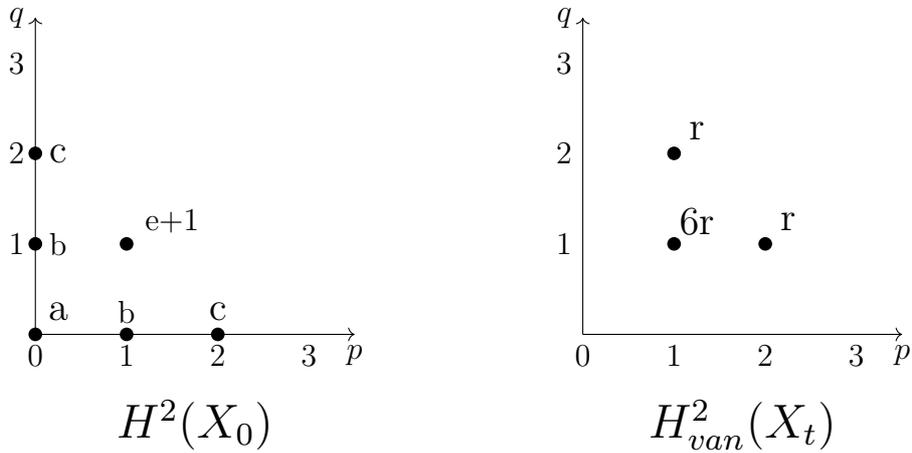
\begin{center}
    \begin{tikzpicture}[scale=.6]
        
            \draw[->] (0,0) -- (7,0) node[anchor=north] {$p$};
            \draw	(0,0) node[anchor=north] {0}
            		(2,0) node[anchor=north] {1}
            		(4,0) node[anchor=north] {2}
            		(6,0) node[anchor=north] {3};
            
            \draw[->] (0,0) -- (0,7) node[anchor=east] {$q$};
            \draw	(0,2) node[anchor=east] {1}
            		(0,4) node[anchor=east] {2}
            		(0,6) node[anchor=east] {3};
            
            \fill (0,0) circle (0.15); 
            \node [scale=1.25] at (.5,.5) {a};
            \fill (2,0) circle (0.15); 
            \node [scale=1.0] at (2,.5) {b};
            \fill (0,2) circle (0.15); 
            \node [scale=1.0] at (.5, 2) {b};
            \fill (4,0) circle (0.15); 
            \node [scale=1.25] at (4,.5) {c};
            \fill (0,4) circle (0.15); 
            \node [scale=1.25] at (.5,4) {c};
            \fill (2,2) circle (0.15); 
            \node [scale=1] at (3,2.5) {e+1};
            
            \node [scale=1.5] at (3.5,-2) {$H^2(X_0)$};
            
            \draw[->] (12,0) -- (19,0) node[anchor=north] {$p$};
            \draw	(12,0) node[anchor=north] {0}
            		(14,0) node[anchor=north] {1}
            		(16,0) node[anchor=north] {2}
            		(18,0) node[anchor=north] {3};
            
            \draw[->] (12,0) -- (12,7) node[anchor=east] {$q$};
            \draw	(12,2) node[anchor=east] {1}
            		(12,4) node[anchor=east] {2}
            		(12,6) node[anchor=east] {3};
            
            \fill (14,2) circle (0.15); 
            \node [scale=1.25] at (14.5,2.5) {6r};
            \fill (16,2) circle (0.15); 
            \node [scale=1.25] at (16.5,2.5) {r};
            \fill (14,4) circle (0.15); 
            \node [scale=1.25] at (14.5,4.5) {r};
            
            \node [scale=1.5] at (15.5,-2) {$H^2_{van}(X_t)$};    
    \end{tikzpicture}
    \end{center}
The exactness of the above sequence, and the fact that $\delta$ has weight $3$ then forces the following two forms of the other two diagrams (in order): 
    
    
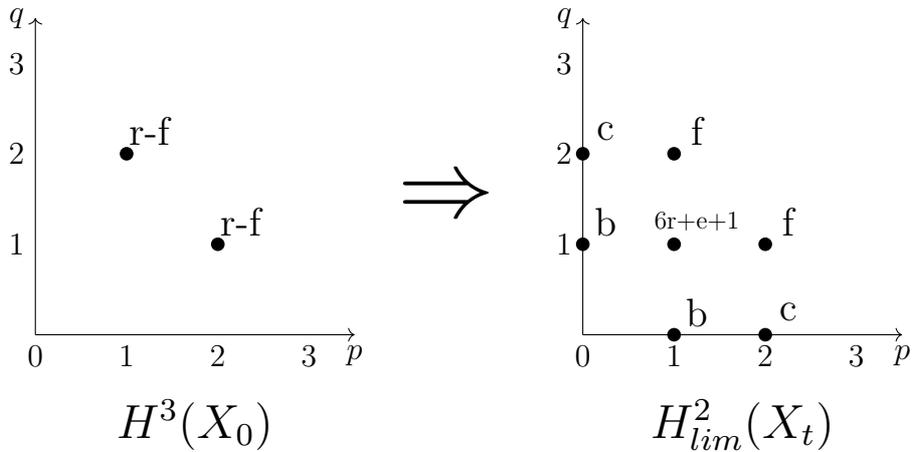
\begin{center}
    \begin{tikzpicture}[scale=.6]
        
            \draw[->] (0,0) -- (7,0) node[anchor=north] {$p$};
            \draw	(0,0) node[anchor=north] {0}
            		(2,0) node[anchor=north] {1}
            		(4,0) node[anchor=north] {2}
            		(6,0) node[anchor=north] {3};
            
            \draw[->] (0,0) -- (0,7) node[anchor=east] {$q$};
            \draw	(0,2) node[anchor=east] {1}
            		(0,4) node[anchor=east] {2}
            		(0,6) node[anchor=east] {3};
            
            \fill (2,4) circle (0.15); 
            \node [scale=1.25] at (2.5,4.5) {r-f};
            \fill (4,2) circle (0.15); 
            \node [scale=1.25] at (4.5,2.5) {r-f};

            \node [scale=1.5] at (3.5,-2) {$H^3(X_0)$};
            
            \node [scale=3] at (9,3) {$\Rightarrow$};
            
            \draw[->] (12,0) -- (19,0) node[anchor=north] {$p$};
            \draw	(12,0) node[anchor=north] {0}
            		(14,0) node[anchor=north] {1}
            		(16,0) node[anchor=north] {2}
            		(18,0) node[anchor=north] {3};
            
            \draw[->] (12,0) -- (12,7) node[anchor=east] {$q$};
            \draw	(12,2) node[anchor=east] {1}
            		(12,4) node[anchor=east] {2}
            		(12,6) node[anchor=east] {3};
            
            \fill (14,2) circle (0.15); 
            \node [scale=.8] at (14.5,2.5) {6r+e+1};
            \fill (16,2) circle (0.15); 
            \node [scale=1.25] at (16.5,2.5) {f};
            \fill (14,4) circle (0.15); 
            \node [scale=1.25] at (14.5,4.5) {f};
            \fill (12,2) circle (0.15); 
            \node [scale=1.25] at (12.5,2.5) {b};
            \fill (12,4) circle (0.15); 
            \node [scale=1.25] at (12.5,4.5) {c};
            \fill (14,0) circle (0.15); 
            \node [scale=1.25] at (14.5,0.5) {b};
            \fill (16,0) circle (0.15); 
            \node [scale=1.25] at (16.5,0.5) {c};

            \node [scale=1.5] at (15.5,-2) {$H^2_{lim}(X_t)$};    
    \end{tikzpicture}
    \end{center}
Now we may deduce from Prop. \ref{prop2.1}[(a),(d)] that:
    
        \begin{align*}
            c + f = h_{3,d}^{2,0} \\
            6r + e + b + f + 1 = h_{3,d}^{1,1}\\
            c + b= h_{3,d}^{0,2} = h_{3,d}^{2,0}
        \end{align*}
So $b=f$ and $6r+2b+e+1 = h_{3,d}^{1,1} \Rightarrow r \leq \frac{1}{6}(h^{1,1}_{3,d} - 1)$
\end{proof}

\begin{prop}\label{prop2.3}
Let $f \in \mathbb{C}[x_0, x_1, \ldots x_n]$ be homogeneous polynomial of degree $d > 2$, let $n=2k+1$, and let $\widetilde{V}(f) \subseteq \mathbb{P}^n$ have only isolated $A_1$ singularities. Then the number $r$ of singular points is bounded by $$ r \leq h_{n,d}^{k,k}-1.$$
\end{prop}

\begin{proof}
     Again we have the following exact sequence of mixed Hodge Structures
    
    \begin{center}
        \begin{tikzcd}[column sep=scriptsize]
            0  \rar &  H^{n-1}(X_0) \rar & H_{lim}^{n-1}(X_t) \rar & H_{van}^{n-1}(X_t) \rar["\delta"] & H^{n}(X_0) \rar & 0
        \end{tikzcd}
    \end{center}
Since $n-1$ is even, we will denote $k=\frac{n-1}{2}$. Since $X_0$ is a compact variety of dimension $n-1$, we know that $H^{n-1}(X_0)$ has the form below by Prop. \ref{prop2.1}(e). Additionally, $H^{2k}_{van}(X_t)$ has Hodge numbers given by $r$ times the Hodge numbers of $H^{n-1}(Y_g)$, where $g: \mathbb{C}^n \rightarrow \mathbb{C}$ has a single $A_1$ singularity. These are calculated explicitly using the formula from \cite{St77} as $h^{k,k}(H^{2k}(Y_g))=1$, and $0$ for the rest. giving the diagram below:
    
    
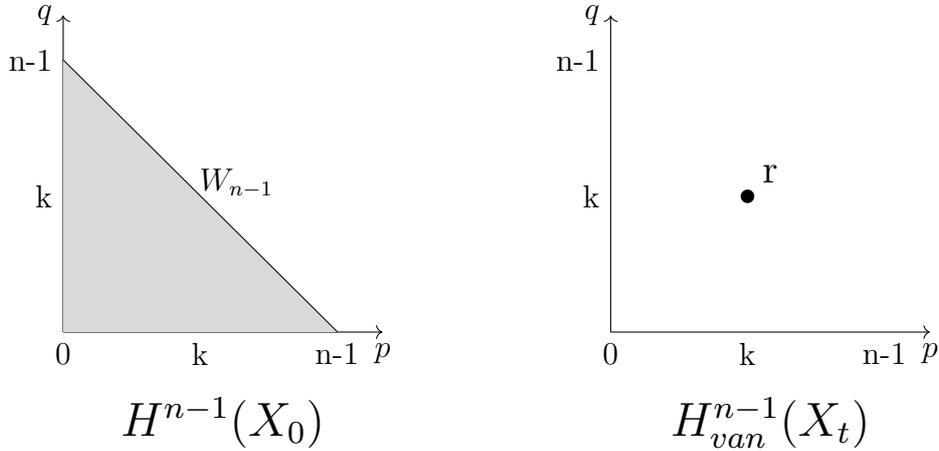
\begin{center}
    \begin{tikzpicture}[scale=.6]
        
             \path[fill=black!15] (0,0) -- (0,6) -- (6,0) -- cycle;
        
            \draw[->] (0,0) -- (7,0) node[anchor=north] {$p$};
            \draw	(0,0) node[anchor=north] {0}
            		(3,0) node[anchor=north] {k}
            		(6,0) node[anchor=north] {n-1};
            
            \draw[->] (0,0) -- (0,7) node[anchor=east] {$q$};
            \draw	(0,3) node[anchor=east] {k}
            		(0,6) node[anchor=east] {n-1};

            \draw[thick]    (0,6) -- (6,0);
            \draw           (3.8,3.3) node {$W_{n-1}$}; 
            \path[fill=black!15] (0,0) -- (0,6) -- (6,0) -- cycle;

            \node [scale=1.5] at (3.5,-2) {$H^{n-1}(X_0)$};

            \draw[->] (12,0) -- (19,0) node[anchor=north] {$p$};
            \draw	(12,0) node[anchor=north] {0}
            		(15,0) node[anchor=north] {k}
            		(18,0) node[anchor=north] {n-1};
            
            \draw[->] (12,0) -- (12,7) node[anchor=east] {$q$};
            \draw	(12,3) node[anchor=east] {k}
            		(12,6) node[anchor=east] {n-1};
            
            \fill (15,3) circle (0.15); 
            \node [scale=1.25] at (15.5,3.5) {r};
            
            \node [scale=1.5] at (15.5,-2) {$H^{n-1}_{van}(X_t)$};    
    \end{tikzpicture}
    \end{center}
The exactness of the above sequence, and the fact that $\delta$ has weight $n$ then forces the following two forms of the other two diagrams (in order): 
    
    
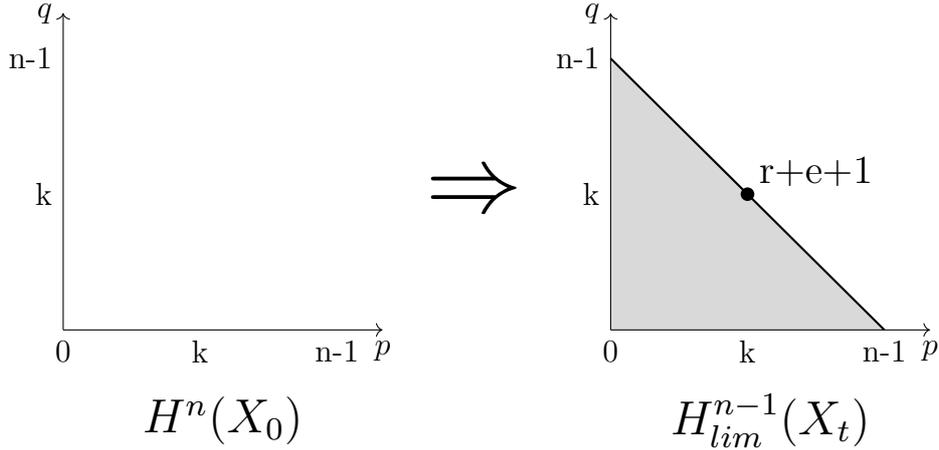
\begin{center}
    \begin{tikzpicture}[scale=.6]
        
            \draw[->] (0,0) -- (7,0) node[anchor=north] {$p$};
            \draw	(0,0) node[anchor=north] {0}
            		(3,0) node[anchor=north] {k}
            		(6,0) node[anchor=north] {n-1};
            
            \draw[->] (0,0) -- (0,7) node[anchor=east] {$q$};
            \draw	(0,3) node[anchor=east] {k}
            		(0,6) node[anchor=east] {n-1};

            \node [scale=1.5] at (3.5,-2) {$H^{n}(X_0)$};
            
            \node [scale=3] at (9,3) {$\Rightarrow$};
            
            \path[fill=black!15] (12,0) -- (12,6) -- (18,0) -- cycle;
            \draw[->] (12,0) -- (19,0) node[anchor=north] {$p$};
            \draw	(12,0) node[anchor=north] {0}
            		(15,0) node[anchor=north] {k}
            		(18,0) node[anchor=north] {n-1};
            
            \draw[->] (12,0) -- (12,7) node[anchor=east] {$q$};
            \draw	(12,3) node[anchor=east] {k}
            		(12,6) node[anchor=east] {n-1};
            

            \draw[thick]    (12,6) -- (18,0);
            \fill (15,3) circle (0.15); 
            \node [scale=1.25] at (16.5,3.5) {r+e+1};

            \node [scale=1.5] at (15.5,-2) {$H^{n-1}_{lim}(X_t)$};    
            
    \end{tikzpicture}
    \end{center}
Now we may deduce from Prop. \ref{prop2.1}[(a),(d)] that:
    
        \begin{align*}
            h^{n-1,0}_{lim} = h_{n,d}^{n-1,0} \\
            r + e + 1 + \sum_{q\neq k}  h^{k,q}_{lim} = h_{n,d}^{k,k}
        \end{align*}
The second equation shows that $r \leq h_{n,d}^{k,k}-1$, as desired. 
\end{proof}

\section{Varchenko's Bound}\label{S3}

Varchenko was the first to integrate the concept of the singularity spectrum in attempt to bound the number of singularities of a projective hypersurface. His original proof can be found in \cite{Var83}, and a further discussion of the proof can be found in \cite{vS20}. The conical bounding method is able to duplicate these bounds in the case of ordinary double points by means of more advanced properties of the spectrum. A discussion of this process is given in the next section. 

We use the convention of the Steenbrink spectrum (denoted $\sigma$, and briefly recalled after the theorem below).  For notational sake, let $\{\sigma\}$ denote the set of spectral summands \emph{with multiplicity}. That is, for the spectrum $\sigma=\left[ \frac{1}{3} \right]+2\left[\frac{1}{2}\right]$ we would have $\{\sigma\}=\left\{\frac{1}{3}, \frac{1}{2}, \frac{1}{2} \right\}$. For any subset $S \subseteq \mathbb{R}$ and spectrum $\sigma$, let $S \cap^\# \{\sigma\} $ count the number of times the sets intersect (for example if $S = \left\{\frac{1}{2}\right\}$ and $\sigma$ is the one above $S \cap^{\#} \{\sigma\} = 2$).  Varchenko's bound can be summarized as follows.

\begin{thm}[\cite{Var83,vS20}]\label{th3.1}
Let $Z \subseteq \mathbb{P}^n$ be a hypersurface of degree $d$, with only isolated singular points $P_1, \ldots, P_r$. Let $g_i: (\mathbb{C}^n,0) \rightarrow (\mathbb{C},0)$, for $1 \leq i \leq r$ denote the corresponding germs defined locally about $P_i$. Then for each $\alpha \in \mathbb{R}$, one has:
    
    $$  (\alpha, \alpha +1)  \cap^{\#} \{\sigma_{x_1 ^d + \ldots + x_n^d,0} \} \geq \sum_{i=1}^r  (\alpha, \alpha +1)  \cap^{\#} \{\sigma_{g_i,0} \}  $$

    \noindent Since $\sigma_{x_1 ^d + \ldots + x_n^d,0} = \gamma_d^{*n}$ where $\gamma_d:= \sum_{i=1}^{d-1}\left[- \frac{i}{d} \right]$ this can be restated as:
    
     $$  (\alpha, \alpha +1)  \cap^{\#} \{ \gamma_{d}^{*n} \} \geq \sum_{i=1}^r  (\alpha, \alpha +1)  \cap^{\#} \{\sigma_{g_i,0} \}  $$
    
\end{thm}

Here we recall that $\sigma_{g,0}=\sum_{q \in \mathbb{Q}} m_{q} [q] \in \mathbb{Z}[\mathbb{Q}]$ means that $$m_{q}=\dim\{ \mathrm{Gr}_F^{\lfloor n-q-1\rfloor}(\tilde{H}^{n-1}(Y_{g,0}))_{e^{-2\pi i q}}\},$$ where the subscript denotes the eigenvalue of the semisimple part $T^{ss}$ of the monodromy operator $T$.   The star notation is defined by $[q]*[q']=[q+q'+1]$ on generators, and the spectrum of $\sum_{k=1}^n x_k^{d_k}$ is given by $\gamma_{d_1}*\cdots *\gamma_{d_n}$.

\begin{exmp}
    We will consider the case when $Z\subseteq \mathbb{P}^n$ is a hypersurface of degree $d$, with only ordinary double points. Then all $g_i$ have the form $g_i=x_1^2+\ldots + x_n^2$ up to analytic equivalence and spectra $\sigma_{g_i,0}= \left[\frac{n}{2}-1\right]$. Therefore any choice of $\alpha \in \left(\frac{n}{2}-2,\frac{n}{2}-1 \right)$ will yield a bound on the number of singularities $r$ possible by the above theorem. 
    
    If $n=3$, then $\sigma_{g_i,0}= \left[\frac{1}{2} \right]$ for all $1\leq i \leq r$.  The Appendix details the spectra $\gamma_d^{*3}$ of $x_1^d+x_2^d+x_3^d$ (note: in the appendix this will correspond to $n=2$). It becomes clear from these spectra that our lowest bound will be obtained by choosing $\alpha=-\frac{1}{2}+ \frac{1}{2d}$ for $d$ even and $\alpha=\frac{\lfloor \frac{d}{2}\rfloor+1}{d}-1$ for general $d$. 
    
    Let $\alpha$ be chosen in this way. Then, in the notation of the theorem, the right hand side of the inequality becomes
    $$ r \left( (\alpha, \alpha + 1 )  \cap^{\#} \{ \sigma_{g_i,0} \} \right) = \# \bigcup_{i=1}^r \left\{ \frac{1}{2} \right\} = r,$$
    where the left hand side will be the number of summands in $\gamma_d^{*3}$ which fall in the interval $(\alpha, \alpha+1))$ inspected from the chart. 
    
    For $d=2$ we get $(\alpha, \alpha+1)=\left( -\frac{1}{4}, \frac{3}{4} \right) \Rightarrow r\leq 1$. For $d=3$ we get $(\alpha, \alpha+1)=\left( -\frac{1}{3}, \frac{2}{3} \right) \Rightarrow r\leq 1+3 =4 $. For $d=7$, $r\leq 6+10+15+21+25+27 = 104 $. In fact it is the case that the bounds match up with those of the next section, at least in the case of $A_1$ singularities.

\end{exmp}

\section{Conical Bounding Method for Pham-Brieskorn}\label{S4}

Throughout this section, $f\in \mathbb{C}[x_0,\ldots,x_n]$ will denote a homogeneous polynomial of degree $d$.  We write as above $\widetilde{V}(f)\subset \mathbb{P}^n$ for the projective hypersurface it defines, and $V(f)\subset \mathbb{C}^{n+1}$ for the affine variety it defines (which is the just the cone on $\widetilde{V}(f)$).  If $\widetilde{V}(f)$ has isolated singularities, then the singularity locus $\Sigma \subset V(f)$ has dimension one and consists of lines through  the origin.  In the neighborhood of an isolated singularity $P\in \widetilde{V}(f)$, we can represent $\widetilde{V}(f)$ in local analytic coordinates (on $\mathbb{P}^n$, about $P$) by the vanishing of a polynomial $g\colon (\mathbb{C}^n,0)\to (\mathbb{C},0)$.  (That is, we are only interested in the intersection of $g=0$ with a small ball about the origin.)

Regarding the definition of the spectrum for a \emph{nonisolated} singularity $p$ (with local affine equation $F=0$ and Milnor fiber $Y$) on a hypersurface of dimension $n$, the $\tilde{H}^k(Y)$ can be nonvanishing for $n-\dim(\Sigma)\leq k\leq n$ (where $\Sigma$ is the singularity locus).  Accordingly, we define the spectrum as an alternating sum $\sigma_{F,p}:=\sum_{j \geq 0} (-1)^{j}\sigma_{F,p}^{n-j}$ where $\sigma^k_{F,p}$ is derived from the MHS and $T^{ss}$-action on $\tilde{H}^k(Y)$, see \cite{St89}.  The main point is that for the cone singularity at the origin of $V(f)$, this takes the form $\sigma_{f,0}=\sigma_{f,0}^n - \sigma_{f,0}^{n-1}$, which may not be effective.  In order to circumvent this problem, we relate $\sigma_{f,0}$ to the (effective) spectra of \emph{isolated} singularities in two different ways, given by the next two theorems.

\begin{thm}[{\cite[Theorem 6.3]{St89}}]\label{th4.1}
Assume that $\widetilde{V}(f) \subseteq \mathbb{P}^n$ has only isolated singularities, $P_1,\ldots P_r$.  Let each germ $g_i: (\mathbb{C}^n,0) \rightarrow (\mathbb{C},0)$ be analytically equivalent to a Pham-Brieskorn polynomial (i.e. a polynomial of the form $\sum_{j=1}^m {x_j}^{a_j}$ , $a_j \in \mathbb{N}$). Let $\mu_i$ denote the Milnor number of $g_i$, and let the values $\lambda_{ij}$ be defined from the spectra by $\sigma_{g_i,0}= \sum_{j=1}^{\mu_i} [\lambda_{ij}]$.  Then there exists a sufficiently general linear form $\ell$ and sufficiently large $k \in \mathbb{N}$ such that $f_k=f+\epsilon \ell^k$ has an isolated singularity at $0$ for $\epsilon\neq 0$ sufficiently small, and
    $$ \sigma_{f_k,0}= \sigma_{f,0} + \sum_{i,j} \bigg[\lambda_{ij}-\frac{\alpha_{ij}}{k} \bigg]* \beta_{k}$$
with $\beta_m= \sum_{i=0}^{m-1} \left[\frac{-i}{m} \right]$ and $\alpha_{ij}=d\lambda_{ij}-\lfloor d \lambda_{ij} \rfloor$.
\end{thm}

The $f_k$ is called a \emph{Yomdin deformation}.  Note that it is \emph{not} necessarily homogeneous.

\begin{thm}[{\cite[Theorem 6.1]{St89}}] \label{th4.2}
With the same notation as in Theorem \ref{th4.1}, but dropping the Pham-Brieskorn assumption on the isolated singularities, we have
        $$ \sigma_{f,0}= \sigma_{h,0} - \sum_{i,j} \bigg[\lambda_{ij}-\frac{\alpha_{ij}}{d} \bigg]* \beta_{d},$$
where $h$ is a homogeneous polynomial of degree $d$ such that $V(h)$ has an isolated singularity at $0$.
\end{thm}

Before continuing we record the following facts:
\begin{lemma}\label{lem4.3}

    \begin{enumerate}[label=(\alph*)]
        \item The Milnor fiber of an $m$-dimensional isolated hypersurface singularity is $(m-1)$-connected, so the spectrum of any germ defined locally about this singularity is effective (i.e. all of its summands' coefficients are nonnegative).
        \item Let $h \in \mathbb{C}[x_0, \ldots , x_n]$ be homogeneous polynomial of degree $d$ and let $h$ have an isolated singularity at $0$. Then $\sigma_{h,0} = \gamma_d^{*(n+1)}$.
    \end{enumerate}

\end{lemma}

We note that in \cite{St89}, it was implicitly assumed that $k=d+1$ is a sufficiently high power of the general linear form $\ell$ in the context of Theorem 4.1. We later prove the more general Lemma 6.2 using the work of \cite{KL1}, verifying this assumption in greater generality.  With this in mind, we arrive at the following bounding argument which expands on an idea of Steenbrink and T. de Jong in the case of $A_1$ singularities.

\begin{thm}[Conical bounding method]\label{thm4.4}
Assume that $\widetilde{V}(f) \subseteq \mathbb{P}^n$ has only $r$ isolated singularities of a single isomorphism class, describable in local coordinates by a Pham-Brieskorn polynomial $g$ with $\sigma_{g,0}=\sum_j [\lambda_j]$.   Define $\alpha_{j}=d\lambda_{j}-\lfloor d \lambda_{j} \rfloor$.  Then we have
     $$ \sigma_{f_{d+1},0}= \gamma_d^{*(n+1)} - r \left( \sum_{j} \bigg[\lambda_{j}-\frac{\alpha_{j}}{d} \bigg]* \beta_{d} - \sum_{j} \bigg[\lambda_{j}-\frac{\alpha_{j}}{d+1} \bigg]* \beta_{d+1} \right),$$
and the effectiveness of the spectrum on the left-hand side bounds the number $r$. 
\end{thm}

\begin{proof}
    
    By Theorem \ref{th4.1}, and the assumption that k=d+1 is sufficient in this case, $f+\epsilon \ell^{d+1}$ has an isolated singularity at the origin, and by Lemma \ref{lem4.3}(a), the spectrum of $f+\epsilon \ell^{d+1}$ at the origin is effective. By Lemma \ref{lem4.3}(b),  Theorem \ref{th4.1}, and Theorem \ref{th4.2} we get 
$$\sigma_{f_{d+1},0}-\textstyle\sum_{i,j}[\lambda_{ij}-\tfrac{\alpha_{ij}}{d+1}]*\beta_{d+1}=\sigma_{f,0}=\gamma_d^{*(n+1)}-\sum_{i,j}[\lambda_{ij}-\tfrac{\alpha_{ij}}{d}]*\beta_d$$
and thus the desired formula for the spectrum of the Yomdin deformation.
\end{proof}

\begin{cor}[given in \cite{St89}]
Assume $\widetilde{V}(f) \subseteq \mathbb{P}^n$ has only $r$ isolated $A_1$ singularities (ordinary double points). Then 
    $$ r \leq  \begin{cases} 
      \text{the coefficient of $\big[\frac{n}{2}-1+\frac{1}{d} \big]$ in $\gamma_d^{*(n+1)}$}, & \text{ $nd$ even} \\
      \text{the coefficient of $\big[\frac{n}{2}-1+\frac{1}{2d} \big]$ in $\gamma_d^{*(n+1)}$}, &\text{ $nd$ odd}. \\
   \end{cases}
    $$

\end{cor}

\begin{proof}
    
    The local normal form of each ordinary double point is given by germs $g_i: (\mathbb{C}^n,0) \rightarrow (\mathbb{C},0)$ given by $z_0^2 +\ldots + z_{n-1}^ 2$. And so $\sigma_{g_i,0}= \sigma_{g_1,0}= \big[\frac{n}{2} - 1 \big]$ for $1 \leq i \leq r$. This yields:
    
     $$ \sigma_{f_{d+1},0}= \gamma_d^{*(n+1)} - r \left(  \bigg[\frac{n}{2} - 1 -\frac{\alpha_{11}}{d} \bigg]* \beta_{d} -  \bigg[\frac{n}{2} - 1 - \frac{\alpha_{11}}{d+1} \bigg]* \beta_{d+1} \right)$$
    
    \noindent If nd is even then $\alpha_{11}=d\big( \frac{n}{2} - 1 \big)-\big\lfloor d\big( \frac{n}{2} - 1 \big) \big\rfloor = 0$ and
    
        $$ \sigma_{f_{d+1},0}= \gamma_d^{*(n+1)} - r \left(  \bigg[\frac{n}{2} - 1  \bigg]* \bigg[ \beta_{d} - \beta_{d+1} \bigg] \right).$$ As one can calculate, the coefficient of $\big[-\frac{d-1}{d}\big] = \big[\frac{1}{d} - 1 \big]$ on the right side side must be $r$. By the effectiveness of the left hand side we have:
        $$r \leq \text{the coefficient of $\bigg[\frac{n}{2}-1+\frac{1}{d} \bigg]$ in $\gamma_d^{*(n+1)}$} $$
    
    \noindent If nd is odd then $\alpha_{11}=d\big( \frac{n}{2} - 1 \big)-\big\lfloor d\big( \frac{n}{2} - 1 \big) \big\rfloor = -\frac{1}{2}$ and
    
         $$ \sigma_{f_{d+1},0}= \gamma_d^{*(n+1)} - r \left(  \bigg[\frac{n}{2} - 1  - \frac{1}{2d} \bigg]* \beta_{d} -  \bigg[\frac{n}{2} - 1 -\frac{1}{2(d+1)}\bigg]* \beta_{d+1}  \right)$$
        The coefficient of $ \bigg[\frac{n}{2} - 1  + \frac{1}{2d} \bigg] $ in $\bigg[\frac{n}{2} - 1  - \frac{1}{2d} \bigg]* \beta_{d}$ is $1$ and the coefficient of $ \bigg[\frac{n}{2} - 1  + \frac{1}{2d} \bigg] $ in $\bigg[\frac{n}{2} - 1  - \frac{1}{2(d+1)} \bigg]* \beta_{d+1}$ is $0$ 
    
    Therefore the coefficient on $ \bigg[\frac{n}{2} - 1  + \frac{1}{2d} \bigg] $ in 
    
    $$r \left(  \bigg[\frac{n}{2} - 1  - \frac{1}{2d} \bigg]* \beta_{d} -  \bigg[\frac{n}{2} - 1 -\frac{1}{2(d+1)}\bigg]* \beta_{d+1}  \right)$$
    must be $r$.
    By the effectiveness of the left hand side we have:
     $$r \leq \text{the coefficient of $\bigg[\frac{n}{2}-1+\frac{1}{2d} \bigg]$ in $\gamma_d^{*(n+1)}$}. $$

\end{proof}

Before continuing, we offer one more application of Theorem \ref{thm4.4}.

\begin{exmp}

    Let $f \in \mathbb{C}[x_0, x_1, x_2, x_3]$ be homogeneous polynomial of degree $d > 3$, and let $\widetilde{V}(f) \subseteq \mathbb{P}^3$ have only isolated $\widetilde{E_6}$ singularities which have normal form $x^3+y^3+z^3$. Then it can be shown that the number $r$ of singular points is bounded by
    $$ r \leq \frac{\text{the coefficient of $\big[1+\frac{1}{d} \big]$ in $\gamma_d^{*4}$}}{7} $$ in a similar fashion to the proof above with the cases $d\equiv p \mod 3$.  
    
    We note that while Theorem \ref{thm4.4} technically implies the following better bound for larger $d$, it cannot be shown without laborious arithmetic, or made apparent without computational facts illustrated later in this paper. If $p = \lfloor \frac{2d}{3} \rfloor + 1$, then
    
    $$ r \leq \frac{\text{the coefficient of $\big[\frac{p}{d} \big]$ in $\gamma_d^{*4}$}}{7}.  $$

\end{exmp}

\section{Eigenvalue Bounding Method}\label{S5}

The statements of Section \ref{S3} might cause one to wonder whether the bound has more to do with the Hodge filtration information in the spectrum of $f_k$ or the multiplicity of eigenvalues of $T^{ss}$, the semisimple part of monodromy.  In this section, we attempt to bound the number of singularities of $\widetilde{V}(f)$ using only the eigenvalue multiplicities in the Milnor fiber cohomology of $f_k$. We find that while this method produces a decent bound for small values of the degree $d$, as $d$ increases the bound on the number of singularities is often not nearly as sharp as the methods of Section \ref{S3}.  However, we already have the characteristic polynomials of the monodromy for $f_k$ for a much larger group of polynomials $f$ provided in \cite{Si90}.

Let $F \in \mathbb{C}[x_0,\ldots, x_n]$ have at most a single isolated singularity at 0. If $Y_F$ denotes the Milnor fiber of $F$, then we know that the reduced homology group $\widetilde{H}_k(Y_F)$ can only be nonzero for $k=n$.  We will denote by $M[F](\lambda)$ the characteristic polynomial of the algebraic monodromy acting on $\widetilde{H}_{n}(Y_F)$. In \cite{Mi68}, Milnor gives $M[F](\lambda)$ for $F=\sum_{i=0}^{n} z_i^d$. Since any homogeneous polynomial of degree $d$ in $n+1$ variables with only an isolated singularity at $0$ is a $\mu$-constant deformation of such a polynomial we have:

\begin{prop}

    Let $F \in \mathbb{C}[x_0, \ldots x_n]$ be homogeneous of degree $d$ with only an isolated singularity at $0$. Then $M[F](\lambda)$ is given by:
    
    $$M[F](\lambda)=  \begin{cases} 
    (\lambda -1)^{-1}(\lambda^d-1)^{\frac{1+(d-1)^{n+1}}{d}} & \text{ $n$ even} \\
     (\lambda -1)(\lambda^d-1)^{\frac{(d-1)^{n+1}-1}{d}} & \text{ $n$ odd} \\
   \end{cases}
    $$
    And since these polynomials are dependent only on $n$ and $d$, we denote this polynomial by $M^{reg}_{n,d}(\lambda)$.

\end{prop}

In \cite{Si90} this polynomial is denoted by $M_d^{reg}(\lambda)$ since $n$ was assumed to be fixed. One may also verify that the spectrum of such a polynomial $\gamma_d^{*(n+1)}$, stated above and in \cite{St89},  is consistent with this characteristic polynomial. 

Before we proceed with a summary of \cite{Si90}, we make a quick note that Siersma insists that his version of a general linear form $\ell$ must be \emph{admissible}, that is $\{\ell=0\} \cap f^{-1}(0)$ has an isolated singularity. However, a close reading of the proof of Lemma \ref{lemKL} detailed later in this paper tells us that our pick of $\ell=y_0$ after a coordinate change will be admissible. 

\begin{lemma}\label{lem5.2}
    Let $f \in \mathbb{C}[x_0, \ldots , x_n]$ be homogeneous polynomial of degree $d>2$ and let $\widetilde{V}(f) \subseteq \mathbb{P}^n$ have only isolated singularities. Then there is a  suitable coordinate transformation such that $f(y_0, \ldots, y_n)$ is a homogeneous polynomial of degree $d$  such that $y_0$ is an admissible linear form.  Furthermore there exists an $\varepsilon >0$ such $f+\varepsilon y_0^{k}$ has an isolated singularity at $\underbar{0}$ for all $k>d$.
\end{lemma}

\begin{proof}
     In the notations of the proof of \ref{lemKL}, $ f^{-1}(0) \cap \text{sing}(f)\cap\{y_0 = 0\}\subseteq  \text{sing}(\pi)\cap\{y_0=0\}={\underbar{0}}$. Therefore $\{y_0=0\} \cap f^{-1}(0)$ can have at most an isolated singularity at $\underbar{0}$. Since $f$ is homogeneous in our case, for $d>2$, $f$ must have a singularity at $\{\underbar{0}\}$. 
\end{proof}

The following is stated in \cite[p195]{Si90}:

\begin{thm}\cite{Si90} \label{th5.3}
Let $f \in \mathbb{C}[x_0, \ldots , x_n]$ be homogeneous polynomial of degree $d$ such that $\widetilde{V}(f) \subseteq \mathbb{P}^n$ has only isolated singularities, $P_1,\ldots P_r$.  Let each germ $g_i: (\mathbb{C}^n,0) \rightarrow (\mathbb{C},0)$ be  defined locally about $P_i$, and let $\mu_i$ and $T_i$ denote its Milnor number and algebraic monodromy operator on $H_{n-1}(X_{g_i})$ respectively. Then for our choices of $\varepsilon$ and an admissible linear form $\ell$, we have

\begin{enumerate}
    \item $M[f+\varepsilon \ell^d](\lambda)=M^{reg}_{n,d}(\lambda)$
    \item For $k > d$, we have
    $$ M[f+\varepsilon \ell^k](\lambda) = \frac{M^{reg}_{n,d}(\lambda)}{(\lambda^d - 1)^{\sum \mu_i}} \cdot \prod_{i=1}^r \det(\lambda^kI-T_i^{k-d} ).$$
\end{enumerate}
\end{thm}

This leads to the following bounding argument which mimics the logic of the conical bounding method above. 

\begin{cor}
With the assumptions and notation of Theorem \ref{th5.3}, the number of singularities is bounded by the following relation:
    
    $$ (\lambda^d - 1)^{\sum \mu_i} \bigg|
    M^{reg}_{n,d}(\lambda) \cdot \prod_{i=1}^r \det(\lambda^{d+1}I-T_i ) $$
  In particular, if each $g_i$ has the same singularity type, then:
    
    $$ (\lambda^d - 1)^{r\mu_1} \bigg|
    M^{reg}_{n,d}(\lambda) \cdot [\det(\lambda^{d+1}I-T_1 )]^r $$
    
\end{cor}

\begin{proof}
$f+\varepsilon \ell^{d+1}$ has an isolated singularity at $\underbar{0}$ by Lemma \ref{lem5.2}. Therefore $M[f+\varepsilon \ell^k](\lambda)$ must not have any poles (as it must be a polynomial). This implies the denominator of (2) in Theorem \ref{th5.3} must divide the numerator. 
\end{proof}

We look at this bound in action with a case we already know from above.

\begin{prop}
 
    Let $f \in \mathbb{C}[x_0, x_1, x_2, x_3]$ be homogeneous polynomial of degree $d > 3$, and let $\widetilde{V}(f) \subseteq \mathbb{P}^3$ have only isolated $\widetilde{E_6}$ singularities. Then the number $r$ of singular points is bounded by
    $$ r \leq \frac{(d-1)^4-1}{8d} = \frac{1}{8}d^3-\frac{1}{2}d^2+\frac{3}{4}d-\frac{1}{2}$$ 

\end{prop}

\begin{proof}

    In this case, $n=3$ is odd, so 
    
    $$ M^{reg}_{3,d} = (\lambda -1)(\lambda^d-1)^{\frac{(d-1)^{4}-1}{d}}$$

    Each $\widetilde{E_6}$ singularity corresponds to $\mu_1=8$, and $\det(\lambda I- T_1) = \frac{(\lambda^3 - 1)^3}{(\lambda - 1)}$ Which can both be inferred from the spectrum. Therefore our bounding argument gives us:
    
    $$ (\lambda^d - 1)^{8r} \bigg| (\lambda -1)(\lambda^d-1)^{\frac{(d-1)^{4}-1}{d}} \cdot \frac{(\lambda^{3(d+1)}-1)^{3r}}{(\lambda-1)^r}  $$
    
    We will count the multiplicity of the eigenvalue $e^{\frac{d-1}{d}\cdot 2\pi i}$ on each side. We note that this is not a root of $\lambda^{3(d+1)}-1)^{3r}$ since this would imply $\frac{3(d+1)(d-1)}{d}$ is an integer $\Rightarrow d | 3(d^2-1) = 3d^2-3 \Rightarrow d | -3$ which is not possible since $d > 3$. Therefore the multiplicity on the left side is $8r$ and the right side is $\frac{(d-1)^4-1}{d}$. This implies $r \leq \frac{(d-1)^4-1}{8d}.$ 
\end{proof}

For $A_1$ singularities, the eigenvalue method gives the following bound:

\begin{prop}
    Let $f \in \mathbb{C}[x_0, x_1, x_2, \ldots x_n]$ be a homogeneous polynomial of degree $d > 2$, and let $\widetilde{V}(f) \subseteq \mathbb{P}^n$ have only isolated $A_1$ singularities. Then the number $r$ of singular points is bounded by
     $$ r \leq 
        \begin{cases}
        \frac{1+(d-1)^{n+1}}{d} & \text{n even} \\
        \frac{(d-1)^{n+1}-1}{d} & \text{n odd}.
        \end{cases}$$ 
\end{prop}

\begin{proof}
Each local $A_1$ singularity corresponds to the spectrum $[\tfrac{n}{2}-1]$, hence to the eigenvalue $(-1)^n$; so the characteristic polynomial of each $T_i$ is $(\lambda -(-1)^n)$.  By our corollary this implies:
\begin{equation*}
(\lambda^d-1)^r\bigl\vert M^{reg}_{n,d}(\lambda)\cdot (\lambda^{d+1}-(-1)^n)
\end{equation*}
Picking out the root $e^{\frac{2\pi i}{d}}$ and counting the multiplicities on each side, we deduce that  $r\leq \tfrac{1}{d}\left((d-1)^{n+1}-(-1)^{n+1}\right)$.
\end{proof}

\section{A generalization of the conical bound}\label{S6}

As we have seen from the last two sections, it is necessary that we consider the spectrum and not just the characteristic polynomial of $T^{ss}$ on $H^n(Y_{f_k})$ to get the sharpest bound on the number of singularities.  This has to do with the fact that the Hodge filtration further sorts the eigenvalues of the monodromy, resulting in smaller multiplicities.  The isolated singularities $P_i$ turn out to contribute to the spectrum of the Yomdin deformation in a subtle way, which involves ``pairing'' the action of both horizontal and vertical monodromies on $H^{n-1}(Y_{g_i})$.  Here ``horizontal'' means that $t$ goes about the origin of the disk, while ``vertical'' means to go about the cone point on the $i$th component of $\Sigma$.

Fortunately, Saito and Siersma have left us with the necessary tools to generalize Theorem \ref{th4.1} in such a way that we can generalize the conical bounding process. We will start off by giving Steenbrink's Conjecture from \cite{St89}, which was later proven by Saito in the vast generality of mixed Hodge modules in \cite{Sa91}, and was later specified in more detail in a context closer to our own in \cite[Thm. 7.5]{KL2}. We further contextualize this theorem in the case of homogeneous polynomials.

\begin{thm}[SS Formula for Homogeneous Cone Case] \label{thm1}

    Let $f \in \mathbb{C}[x_0, \ldots , x_n]$ be homogeneous polynomial of degree $d$ such that $\widetilde{V}(f) \subseteq \mathbb{P}^n$ has only isolated singularities, $P_1,\ldots P_r$. Let each germ $g_i: (\mathbb{C}^n,0) \rightarrow (\mathbb{C},0)$ be defined locally about $P_i$, $\mu_i$ denote the Milnor number of $g_i$, and write $\sigma_{g_i,0}= \sum_{j=1}^{\mu_i} [\lambda_{ij}]$. Put $\alpha_{ij}= d\lambda_{ij}-\lfloor d\lambda_{ij} \rfloor$. Then for a sufficiently general linear form $\ell$ and small $\varepsilon\neq 0$, $f_k=f+\varepsilon \ell^k$ has an isolated singularity at $0$ and
     $$ \sigma_{f_{k},0}= \gamma_d^{*(n+1)} - \sum_{i,j} \bigg[\lambda_{ij}-\frac{\alpha_{ij}}{d} \bigg]* \beta_{d} + \sum_{i,j} \bigg[\lambda_{ij}-\frac{\alpha_{ij}}{k} \bigg]* \beta_{k}$$
for any $k>d$.    
\end{thm}
The proof is given later in this section.

\vspace{3mm}

Let $U\subseteq \mathbb{C}^{n+1}$ be a ball about $0$ and $F\colon U \to \Delta$ a holomorphic germ, with $V_t:=F^{-1}(t)$ smooth for $t\neq 0$.  Assume that $Z:=\mathrm{sing}(V_0)$ has dimension $1$, and let $Z_i$ be its irreducible components.  We assume that their only intersection point is the origin.


\begin{lemma}[{\cite[\S7.2]{KL2}}]\label{lemKL}

Consider a (sufficiently general) linear form $\ell$ on $\mathbb{C}^{n+1}$ such that $\{\ell=0\}\cap Z_i$ is finite for each $i$.  Let $\pi=(F,\ell):U \rightarrow \Delta^2$, and denote the Jacobian matrix by $d\pi$.  Let $\mathcal{P} \subseteq U$ denote the set at which $d\pi$ has rank 1; that is, $\mathcal{P}$ is the intersection of the zero loci of all of the $2 \times 2$ minor determinants. Assume in addition that $\mathcal{P} \cap \pi^{-1}(\Delta\times \{(s=)\,0\}) = \{\underbar{0}\}$.  Writing $\pi_s\colon \mathcal{P}\cap \ell^{-1}(s)\to \Delta \times \{s\}$ for the restriction, every irreducible curve $C_q \subset \cup_{s \in \Delta^*_s}\emph{im}(\pi_s)$ has a parametrization of the form $(t_{C_q}(s),s)$, where (for some $r_{C_q}\in \mathbb{Q}_{\geq 0}$ and $\gamma_{C_q}\in \mathbb{C}^*$)
        $$t_{C_q}(s)=\gamma_{C_q}s^{r_{C_q}} + \emph{higher order terms}.$$
Furthermore if
    $$ r := \max_q\{r_{C_q},0\} \in \mathbb{Q}_{\geq0}, $$
then for every $a>r$,  $F+\ell^{a}$ has an isolated singularity at $0$ (including the vacuous case where it is nonsingular at and near $0$).
    
\end{lemma}

We make a quick note that in the previous lemma, $Z \subseteq \cup_s (\mathcal{P}\cap \ell^{-1}(s))$ but the reverse inclusion need not hold.

\begin{exmp}
    Consider the polynomial $F=x_0(x_1^2+x_2^2) \in \mathbb{C}[x_0, x_1, x_2]$. Then $F$ factors as $x_0(x_1+ix_2)(x_1-ix_2)$, and the set $Z=\text{sing}(X_0)= \{(0,z,-iz)\} \cup \{(0,z,iz)\} \cup \{(z,0,0)\}$ for $z \in \mathbb{C}$. We choose $U$ to be a ball about $\underline{0}$, and see that in $U$, $Z$ has dimension $1$. We see that $Z=Z_1 \cup Z_2 \cup Z_3$ on $U$ where $Z_1=\{(0,z,-iz)\}, Z_2=\{(0,z,iz)\}, Z_3 = \{(z,0,0)\}$. Of course, each $Z_i$ is irreducible of dimension 1. 
    
    We choose $\ell=x_0+x_1$, and note $\{\ell=0\} \cap Z_i = \{\underline{0}\}$ for $i=1,2,3$. With $\pi$ as above, we get:
    $$ d\pi = \begin{bmatrix}
    x_1^2+x_2^2 & 2x_0x_1 &  2x_0x_2\\
    1 & 1 & 0
    \end{bmatrix}$$
    Letting $F_{ij}= | \text{Col i} : \text{Col j} |$ for  $ i < j $, we have  $d\pi$ has rank 1 on $\mathcal{P}=\{F_{12}=0\}\cap \{F_{13}=0\} \cap \{F_{23}=0\}$. Here $\mathcal{P}= Z \cup \{z,2z,0\}$. We verify $\mathcal{P} \cap \pi^{-1}(\Delta_t\times \{0\}) = \{Z \cup \{z,2z,0\} \} \cap \{(y_1,-y_1,y_2)\}=\{\underline{0}\}$ where $y_i \in \mathbb{C}$. We have $\mathcal{P} \cap \ell^{-1}(s)= \{ (\frac{s}{3},\frac{2s}{3},0) \}$, and so $\text{im}(\pi_s)= \{(\frac{4s^3}{27}),s \}$. Therefore our only irreducible curve $C\subset \cup_{s \in \Delta^*_s}\text{Im}(\pi_s)$ is given by $t_C=\frac{4}{27}s^3$, and so $r=3$. It follows that for every $a>3$, $F+\ell^a = x_0(x_1^2+x_2^2) + (x_0+x_1)^a $ has an isolated singularity at $\underline{0}$.
\end{exmp}

\begin{exmp}
    An example of the vacuous case is given by $F=x_0(x_1-x_2^2) \in \mathbb{C}[x_0, x_1, x_2]$. Here, $Z=\text{sing}(V_0)= \{(0,z^2,z)\}$ for $z \in \mathbb{C}$. Again we choose $U$ to be a ball about $\underline{0}$, and here $Z$ itself is irreducible of dimension 1. 
    
    We choose $\ell=x_2$, and note $\{\ell=0\} \cap Z = \{\underline{0}\}$.  With $\pi$ as above, we get
    $$ d\pi = \begin{bmatrix}
    x_1-x_2^2 & x_0 &  -2x_0x_2\\
    0 & 0 & 1
    \end{bmatrix},$$
which has rank 1 on $\mathcal{P}= Z$. We verify $\mathcal{P} \cap \pi^{-1}(\Delta_t,\{0\}) = Z \cap \{(y_1,0,0)\}=\{\underline{0}\}$ where $y_i \in \mathbb{C}$. We have $\mathcal{P} \cap \ell^{-1}(s)= \{ (0,s^2,s) \}$, and so $\text{Im}{\pi_s}= \{0,s \}$. Therefore there is no irreducible curve $C\subset \cup_{s \in \Delta^*_s}\emph{Im}(\pi_s)$ and so $r=0$. We conclude that for every $a>0$, $F+\ell^a = x_0(x_1-x_2^2) + (x_2)^a $ has an isolated singularity at $\underline{0}$. In this case, $x_0(x_1-x_2^2) + (x_2)^a$ is in fact nonsingular for $a>0$. 
    
\end{exmp}

We now apply Lemma \ref{lemKL} to our case of interest.

\begin{lemma}\label{lem6.5}

     Let $f \in \mathbb{C}[x_0, \ldots , x_n]$ be homogeneous polynomial of degree $d$ and let $\widetilde{V}(f) \subseteq \mathbb{P}^n$ have only isolated singularities. Then there exists a sufficiently general linear form $\ell$ such that $f+\ell^{k}$ has an isolated singularity at $\underbar{0}$ for all $k>d$. 

\end{lemma}

\begin{proof}
    
    Recall that the singular locus $\Sigma:=\mathrm{sing}(V(f))$ of the \emph{affine cone} is a union of lines passing through the origin.  Then there exists a suitable change of coordinates $(y_0,\ldots, y_n)$ and a ball $B=B_\epsilon(\underbar{0})$ such that the components $\Sigma_i$ of $\Sigma$ are parametrized by $y_0$. That is, each $\Sigma_i$ has the form:
    
    $$\Sigma_i=\cup_{s\in\Delta_s} \{(s,f^j_i(s),\ldots,f^n_i(s) )\mid f_i^j(0)=0 \quad \forall i \}$$
    Furthermore we can set $(y_1,\ldots, y_n)=\underbar{y}$ and rewrite $f(y_0,\underbar{y})=\sum_{j=0}^d g_j(\underbar{y})y_0^{d-j}$, and since $\widetilde{V}(f)$ has only isolated singularities, we may further assume our coordinates were chosen so that $f(0,\underbar{y})=g_d(\underbar{y})$ defines a smooth hypersurface $\widetilde{V}(g_d)\subseteq \mathbb{P}_{(y_0=0)}^{n-1}$.
    We choose $\ell=y_0$ and let $\pi=(f,\ell): B \rightarrow \Delta^2_{t,s} $. Denote $\partial_if=\partial_{y_i}f$. Then we have:
    $$ d\pi = \begin{bmatrix}
    \partial_0f & \partial_2 f & \ldots & \partial_n f\\
    1 & 0 & \ldots & 0
    \end{bmatrix},$$
    which has rank 1 precisely when $\partial_if=0$ for every $i>0$. 
    
    We must show that $\text{sing}(\pi)\cap\{y_0=0\}=\{\underbar{0}\}$. Let $p=(p_0,\underline{p})\in \text{sing}(\pi)\cap\{y_0=0\}$. Then $p=(0,\underline{p})$ and 
    
    $$\partial_if(p)=\partial_i\Bigg[\sum_{j=0}^d g_j(\underbar{x})x_0^{d-j}\Bigg](p)=0 \text{  for  } i>0.$$ 
    This yields $\partial_ig_d(\underline{p})=0$ for $i>0$ Since $\widetilde{V}(g_d)$ is smooth in $\underbar{x}$, this implies $\underline{p}=\underline{\underline{0}} \Rightarrow p=\underbar{0}.$ Therefore there will exist a $k$ such that for $\varepsilon$ small enough $f+\varepsilon y_0^k$ will have an isolated singularity at $\underbar{0}$.

    In the notation above, $\mathcal{P}=\{\partial_if=0 \mid i>0 \} \cap B$. By Euler's identity, we have $d\cdot f= \sum_{i=0}^n y_i \partial_if$. Any curve $C$ given by the image of $\pi_s$ of points $p_s \in \mathcal{P} \cap \{y_0=s\} \cap \{f=0 \}$ is just given by $t_C=0$. So we consider the case when 
    $d\cdot f = x_0\partial_0f $. The points where $f=0$ just give the curve $t_{C_q}=0$. We let $p_s \in \mathcal{P} \cap \{f \neq 0\}$. That is, $p_s=(s,\underline{p_s})$. By Euler's homogeneous function theorem, 
    $$f(p_s)= \frac{1}{d} \sum_{i=0}^n {p_s}_i \partial_i f (p_s) f(p_s)= \frac{s}{d} \partial_0 f (p_s) = \frac{s}{d} \sum_{j=0}^{d-1} (d-j)g_{j}(\underline{p_s})s^{d-j-1}= \frac{1}{d} \cdot \frac{\text{d}}{\text{ds}} f(s,\underline{p_s})$$
    Let $h(s)=f(s,\underline{p_S})$. Then by above, $h(s)=\frac{d}{s} \cdot h'(s)$ for $s \neq 0$. Therefore:
    $$ \int \frac{h'(s)}{h(s)} ds = d \int \frac{ds}{s} \Rightarrow \log(h(s))=d\log(s) + C \Rightarrow h(s)=As^d$$
    And so, by Lemma \ref{lemKL}, we have shown that $f+\varepsilon y_0^{d+1}$ must have an isolated singularity at $\underbar{0}$.
\end{proof}

The following is detailed in \cite[p195]{Si90}:

\begin{lemma}\label{lemS}
     Let $f \in \mathbb{C}[x_0, \ldots , x_n]$ be homogeneous polynomial of degree $d$ such that $\widetilde{V}(f) \subseteq \mathbb{P}^n$ has only isolated singularities, $P_1,\ldots P_r$.  Let each germ $g_i: (\mathbb{C}^n,0) \rightarrow (\mathbb{C},0)$ be  defined locally about $P_i$. Let $\mu_i$ denote the Milnor number of $g_i$. Let $T_i$ and $\tau_i$ be the algebraic horizontal and vertical monodromy operators respectively, each corresponding to $g_i$. Then we must have:
     $$ {T_i}^{-d} = \tau_i.  $$
\end{lemma}

\begin{proof}[Proof of Theorem \ref{thm1}]
    We note that our choice of $\ell$ and $k>d$ yields an $f_k$ with isolated singularity by Lemma \ref{lem6.5}. Additionally, this lemma was proven in accordance with Lemma \ref{lemKL}, which is precisely the condition contained in the preface to \cite[Thm. 7.5]{KL2}. Therefore our choice of $d$ as the bounding exponent $r$ is sufficient to invoke Saito-Steenbrink, but where the $\alpha_{ij}$ are \emph{a priori} given by the eigenvalues of the vertical monodromy operators.
    
    Recall that there exist matrix representations of $T_i$ and $\tau_i$ in Jordan normal form, so that they each have a Jordan-Chavalley decomposition into the product of a unipotent and semisimple matrices:
    \begin{align*}
        T_i = T_i^{ss}T_i^{u} \\
        \tau_i= \tau_i^{ss}\tau_i^{u}
    \end{align*}
    Furthermore, these representations can be chosen in such a way that there exists a simultaneous eigenbasis $v_{ij}$ for $T_i^{ss}$ and $\tau_i^{ss}$ for which:

    \[       
    \begin{array}{lcl} 
        { T_i}^{ss}v_{ij}=(e^{-2\pi i \lambda_{ij}})v_{ij} & &  \\
        {\tau_i}^{ss}v_{ij}=(e^{2\pi i \alpha_{ij}})v_{ij}  & \mbox{for} & \alpha_{ij} \in [0,1)
    \end{array}
    \]
     Note that while the values $\lambda_{ij}$ are the spectral summands above, this relation is what defines the values $\alpha_{ij}$. By Lemma \ref{lemS}, $T_i^{-d}=\tau_i \Rightarrow (T_i^{ss}T_i^{u})^{-d}=\tau_i^{ss}\tau_i^{u}$, but since the pieces of the Jordan-Chavalley decomposition commute, the LHS is just $(T_i^{ss})^{-d}(T_i^{u})^{-d}$. By the uniqueness of the decomposition into semisimple and unipotent parts, we have $(T_i^{ss})^{-d}=\tau_i^{ss}$. This implies our values of $\alpha_{ij}=d\lambda_{ij}-\lfloor d\lambda_{ij} \rfloor$. 

Finally, by \cite[Thm. 7.5]{KL2}, we must have:
    $$ \sigma_{f_k,0}= \sigma_{f,0} + \sum_{i,j} \bigg[\lambda_{ij}-\frac{\alpha_{ij}}{k} \bigg]* \beta_{k}$$
where $\beta_m= \sum_{i=0}^{m-1} \bigg[- \frac{i}{m} \bigg]$.
Combining this with Theorem \ref{th4.2} now gives the desired formula:
     $$ \sigma_{f_{k},0}= \gamma_d^{*(n+1)} - \sum_{i,j} \bigg[\lambda_{ij}-\frac{\alpha_{ij}}{d} \bigg]* \beta_{d} + \sum_{i,j} \bigg[\lambda_{ij}-\frac{\alpha_{ij}}{k} \bigg]* \beta_{k}.$$
\end{proof}

Considering the case where $k=d+1$ we get the following general bound for multiple singularity types at once: 

\begin{thm}[Generalized Conical bounding method]\label{th6.7}

    Let Let $f \in \mathbb{C}[x_0, \ldots , x_n]$ be a homogeneous polynomial of degree $d$, and let $\widetilde{V}(f) \subseteq \mathbb{P}^n$ have only $r$ isolated singularities given locally by  $g_i:(\mathbb{C}^n,0)\rightarrow (\mathbb{C},0)$ having corresponding Milnor numbers $\mu_i$ for $1\leq i \leq r$. Let $\sigma_{g_i,0}= \sum_{j=1}^{\mu_i} [\lambda_{ij}]$ be their corresponding spectra. Then:
   $$ \sigma_{f_{d+1},0} = \gamma_d^{*(n+1)} - \left( \sum_{i=1}^r\sum_{j=1}^{\mu_i}\sum_{k=1}^{d} \bigg[\frac{\lfloor d\lambda_{ij} \rfloor + k}{d} \bigg] - \sum_{i=1}^r\sum_{j=1}^{\mu_i}\sum_{k=1}^{d+1} \bigg[\frac{\lambda_{ij}+\lfloor d\lambda_{ij} \rfloor + k  }{d+1} \bigg] \right) $$
    and the effectiveness of this spectrum restricts the set of $r$ singularities which can be present.  

\end{thm}

\begin{proof}

We know from Theorem \ref{thm1} and Lemma \ref{lem4.3} that $k=d+1$ gives an effective spectrum which satisfies: 

 $$ \sigma_{f_{d+1},0}= \gamma_d^{*(n+1)} - \sum_{i,j} \bigg[\lambda_{ij}-\frac{\alpha_{ij}}{d} \bigg]* \beta_{d} + \sum_{i,j} \bigg[\lambda_{ij}-\frac{\alpha_{ij}}{d+1} \bigg]* \beta_{d+1}.$$
Since $\alpha_{ij}=d\lambda_{ij}-\lfloor d \lambda_{ij} \rfloor$, this is just

 $$= \gamma_d^{*(n+1)} - \left( \sum_{i,j}\sum_{k=1}^{d} \bigg[\frac{\lfloor d\lambda_{ij} \rfloor + k}{d} \bigg] - \sum_{i,j}\sum_{k=1}^{d+1} \bigg[\frac{\lambda_{ij}+\lfloor d\lambda_{ij} \rfloor +k  }{d+1} \bigg] \right).  $$
\end{proof}

\section{Some more user-friendly formulas}\label{S7}

In this section, we will explain how to adapt Theorem \ref{th6.7} to a formula which serves the purpose of reducing the number of calculations. The caveat is that the formula holds no deeper meaning within greater spectral theory. As one can see, Theorem \ref{th6.7} explicitly describes a relationship between the spectrum $\sigma_{f_{d+1},0}$ and the spectra of local singularities. If we throw away the concept of spectra all together, we are left with theorems which only describe relationships of elements in $\mathbb{Z}{[\mathbb{Q}]}$. 

This becomes more convenient because it allows us to simply throw things away that don't matter to the arithmetic we need to do to simply bound the possible singularities. We get the following:

\begin{thm}\label{th7.1}
    Let Let $f \in \mathbb{C}[x_0, \ldots , x_n]$ be a homogeneous polynomial of degree $d$, and let $\widetilde{V}(f) \subseteq \mathbb{P}^n$ have only $r$ isolated singularities given locally by  $g_i:(\mathbb{C}^n,0)\rightarrow (\mathbb{C},0)$ having corresponding Milnor numbers $\mu_i$ for $1\leq i \leq r$. Let $\sigma_{g_i,0}= \sum_{j=1}^{\mu_i} [\lambda_{ij}]$ be their corresponding spectra. Then the following sum in $\mathbb{Z}{[\mathbb{Q}]}$ is effective:
   $$ \gamma_d^{*(n+1)} - \left( \sum_{i=1}^r\sum_{j=1}^{\mu_i}\sum_{k=1}^{d-1} \bigg[\frac{\lfloor d\lambda_{ij} \rfloor + k}{d} \bigg] + \sum_{\substack{i=1 \\ \\ d\lambda_{ij}}}^r\sum_{\substack{j=1 \\ \\ \notin \mathbb{Z}}}^{\mu_i} \bigg[\frac{\lfloor d\lambda_{ij} \rfloor  }{d} +1 \bigg] \right) $$
    and the effectiveness of this sum restricts the set of $r$ singularities which can be present.  
\end{thm}

We will first need to prove a very short lemma:

\begin{lemma}
    If the quantities:
    $$ \frac{\lfloor d\lambda_{i,j} \rfloor +k }{d} =\frac{\lfloor d\lambda_{gh} \rfloor + \ell + \lambda_{gh}}{d+1}$$ for values $d \in \mathbb{N}$ , $1 \leq \ell \leq d+1$, and $1\leq k \leq d$, then we must have the following:
    \begin{enumerate}
        \item $\ell=d+1$
        \item $d\lambda_{gh} \in \mathbb{Z}$
    \end{enumerate}
\end{lemma}

\begin{proof}
    Assume:
    $$ \frac{\lfloor d\lambda_{i,j} \rfloor +k }{d} =\frac{\lfloor d \lambda_{gh} \rfloor + \ell + \lambda_{gh}}{d+1}$$ 
    $\Rightarrow (d+1)(\lfloor d\lambda_{ij}\rfloor +k)=d(\lfloor d \lambda_{gh} \rfloor + \ell + \lambda_{gh})$. Since the left is an integer, the right must be $\Rightarrow d\lfloor d \lambda_{gh} \rfloor + d\ell + d\lambda_{gh} \in \mathbb{Z} \Rightarrow d\lambda_{gh} \in \mathbb{Z} \Rightarrow \lfloor d\lambda_{gh}\rfloor = d\lambda_{gh}$. This restricts the equality to be:
    
    $\Rightarrow (d+1)(\lfloor d\lambda_{ij}\rfloor +k)=d( d\lambda_{gh}  + \ell + \lambda_{gh})= d(d+1)\lambda_{gh}+d\ell$
    
    $\Rightarrow \lfloor d\lambda_{ij}\rfloor +k= d\lambda_{gh}+\frac{d}{d+1}\ell$ Where the left must be an integer so the right must be. Since $d\lambda_{gh}$ is also an integer
    
    $\Rightarrow \frac{d}{d+1} \ell \in \mathbb{Z}$. Our values of $\ell$ only range $1\leq \ell \leq d+1 \Rightarrow \ell=d+1$.
\end{proof}

\begin{proof}[Proof of Theorem 7.1]

    By Theorem \ref{th6.7}, we have that
    $$ \sigma_{f_{d+1},0} = \gamma_d^{*(n+1)} - \left( \sum_{i=1}^r\sum_{j=1}^{\mu_i}\sum_{k=1}^{d} \bigg[\frac{\lfloor d\lambda_{ij} \rfloor + k}{d} \bigg] - \sum_{i=1}^r\sum_{j=1}^{\mu_i}\sum_{k=1}^{d+1} \bigg[\frac{\lambda_{ij}+\lfloor d\lambda_{ij} \rfloor + k  }{d+1} \bigg] \right) $$
    is an effective sum in $\mathbb{Z}{[\mathbb{Q}]}$. By our lemma, the only summands of the right triple sum which may cancel with those of the left triple sum in the subtraction are those that satisfy the properties in the conclusion of the lemma. Therefore
    
    $$\gamma_d^{*(n+1)} + \sum_{\substack{i=1 \\ \\ d\lambda_{ij}}}^r\sum_{\substack{j=1\\ \\ \in \mathbb{Z}}}^{\mu_i} \bigg[\frac{\lambda_{ij}+ d\lambda_{ij}}{d+1}   +1 \bigg]- \sum_{i=1}^r\sum_{j=1}^{\mu_i}\sum_{k=1}^{d} \bigg[\frac{\lfloor d\lambda_{ij} \rfloor + k}{d} \bigg]  $$
    $$=\gamma_d^{*(n+1)} + \sum_{\substack{i=1 \\ \\ d\lambda_{ij}}}^r\sum_{\substack{j=1\\ \\ \in \mathbb{Z}}}^{\mu_i} \bigg[\lambda_{ij} +1 \bigg]- \sum_{i=1}^r\sum_{j=1}^{\mu_i}\sum_{k=1}^{d} \bigg[\frac{\lfloor d\lambda_{ij} \rfloor + k}{d} \bigg]  $$
is an effective sum in $\mathbb{Z}[\mathbb{Q}]$.
But

$$ \sum_{i=1}^r\sum_{j=1}^{\mu_i}\sum_{k=1}^{d} \bigg[\frac{\lfloor d\lambda_{ij} \rfloor + k}{d} \bigg] =  \sum_{\substack{i=1 \\ \\ d\lambda_{ij}}}^r\sum_{\substack{j=1\\ \\ \in \mathbb{Z}}}^{\mu_i} \bigg[\lambda_{ij}+ 1 \bigg] + \sum_{\substack{i=1 \\ \\ d\lambda_{ij}}}^r\sum_{\substack{j=1\\ \\ \notin \mathbb{Z}}}^{\mu_i} \bigg[\lambda_{ij}+ 1 \bigg]+ \sum_{i=1}^r\sum_{j=1}^{\mu_i}\sum_{k=1}^{d-1} \bigg[\frac{\lfloor d\lambda_{ij} \rfloor + k}{d} \bigg] $$
So we conclude that
$$ \gamma_d^{*(n+1)} - \left( \sum_{i=1}^r\sum_{j=1}^{\mu_i}\sum_{k=1}^{d-1} \bigg[\frac{\lfloor d\lambda_{ij} \rfloor + k}{d} \bigg] + \sum_{\substack{i=1 \\ \\ d\lambda_{ij}}}^r\sum_{\substack{j=1 \\ \\ \notin \mathbb{Z}}}^{\mu_i} \bigg[\frac{\lfloor d\lambda_{ij} \rfloor  }{d} +1 \bigg] \right) $$
is effective in $\mathbb{Z}{[\mathbb{Q}]}$.
\end{proof}

The power of this form, as opposed to the one in Theorem \ref{th6.7} is immense. The original theorem would have one believe, on first glance, that it were possible to have a set of two types of local singularities $g_1$ and $g_2$ such that the values of the summands $ \left[ \frac{ \lambda_{1j} + \lfloor d\lambda_{1j} \rfloor +k}{d+1} \right]$ with positive coefficients corresponding to $g_1$ give extra wiggle room to the coefficients of $\gamma_d^{*(n+1)}$ to cancel out the summands  $\left[\frac{\lfloor d\lambda_{2j} \rfloor + k}{d} \right]$ with negative coefficients corresponding to $g_2$. With our theorem in this section, we have proven that this possibility is, in fact, irrelevant to our use of the bounding method. 

We prove the following statements, which serve to further bridge the gap of how our bound becomes increasingly similar to that of Varchenko in Section \ref{S3}. 

\begin{lemma}

    The number of positive integer solutions $(x_1,\ldots , x_k)$ to the equation:
        $$ \sum_{i=1}^k x_i = N $$
    for some positive integer $N$ subject to the constraints $1 \leq x_i \leq \alpha$ for $i=1, \ldots, k$ is given by:
        $$ \sum_{i=0}^m (-1)^i \binom{k}{i} \binom{N-\alpha i -1}{k-1} $$
    where $m = \Big\lfloor \min \big\{ k, \frac{N-k}{\alpha} \big\} \Big\rfloor$ acts as a truncator. 
    
    Furthermore, in the case that $\min \big\{ k, \frac{N-k}{\alpha} \big\} = \frac{N-k}{\alpha}$ we may choose any $  \left\lfloor \frac{N-k}{\alpha} \right\rfloor \leq m \leq \left\lfloor \frac{N-1}{\alpha} \right\rfloor $, as this simply adds zero terms in the sum.   
\end{lemma}

\begin{proof}
    This follows from a basic combinatorial argument using a ``stars and bars'' style proof and inclusion exclusion principles. 
\end{proof}

\begin{prop}\label{prop7.4}
    For any $p \in \mathbb{Z}$ such that $n-d \leq p \leq n(d-1)-d$, we have:
    $$ \left\{ \frac{p}{d} \right\} \cap^{\#} \left\{\gamma_d^{*n} \right\} = \sum_{i=0}^{\left\lfloor n-1 - \frac{p}{d-1} \right\rfloor}(-1)^i\binom{n}{i} \binom{d(n-1)-p-1 -(d-1)i}{n-1}$$
    and this completely determines the spectrum $\gamma_d^{*n}$.
\end{prop}

\begin{proof}
    
    Recall that  $$ \gamma_d^{*n}= \left( \sum_{i=1}^{d-1} \left[ -\frac{i}{d} \right] \right)^{*n} =  \sum_{x_1,\ldots,x_{n}=1}^{d-1} \left[ n -1 - \frac{\sum_{i=1}^{n} x_i}{d} \right]$$
    
    So it's clear that every summand of $\gamma_d^{*n}$ has the form $\frac{p}{d}$ where the numerator $p \in \mathbb{Z}$ is restricted to the range $n-d \leq p \leq n(d-1)-d$. Calculating the coefficient of each summand amounts to counting the number of ways the sum $\sum_{i=1}^{n} x_i = d(n-1)-p$ subject to the constraint $1 \leq x_i \leq d-1 $ for $1\leq i \leq n$. The result then immediately follows from the above lemma.
\end{proof}

We make note of the following cute fact, which follows immediately from Theorem \ref{th1.1} and Proposition \ref{prop7.4}:

\begin{cor}

    Let $d>n$ and the values $\left[h_{n,d}^{k, n-1-k} \right]'$ be the primitive hodge numbers of a smooth hypersurface in $\mathbb{P}^n$ of degree $d$. Then for $k \leq \frac{n-1}{2}$, we have:
    
    $$ \left[h_{n,d}^{k, n-1-k} \right]' = \left\{ n-k-1 \right\} \cap^{\#} \{ \gamma_d^{*(n+1)}\}. $$

\end{cor}

We state one more lemma comparing the coefficients of $\gamma_d^{*n}$ and $\gamma_d^{*(n+1)}$:

\begin{lemma}
    The following is a result of combinatorial arithmetic:
    $$  \left\{ \frac{p}{d}  \right\}  \cap^{\#} \{ \gamma_d^{*(n+1)} \}= \left\{ \frac{p-1}{d}, \frac{p-2}{d}, \ldots , \frac{p-(d-1)}{d} \right\} \cap^{\#} \{ \gamma_d^{*n} \} $$
    $$ =  \left(\frac{p}{d}-1, \frac{p}{d} \right)  \cap^{\#} \{ \gamma_d^{*n} \} $$
\end{lemma}

The following is an equivalent statement of the conical bound:

\begin{thm}[Alternative statement of the conical bound] \label{th7.7}
    Let $f \in \mathbb{C}[x_0, \ldots , x_n]$ be a homogeneous polynomial of degree $d$, and let $\widetilde{V}(f) \subseteq \mathbb{P}^n$ have only $r$ isolated singularities given locally by  $g_i:(\mathbb{C}^n,0)\rightarrow (\mathbb{C},0)$ having corresponding Milnor numbers $\mu_i$ for $1\leq i \leq r$. Let $\sigma_{g_i,0}$ be their corresponding spectra. Then for every $p \in \mathbb{Z}$ we must have:
      $$\left\{ \frac{p}{d}  \right\} \cap^{\#} \{\gamma_d^{*(n+1)} \}  \geq \sum_{i=1}^r  \left( \frac{p}{d}-1 , \frac{p}{d}\right) \cap^{\#} \{ \sigma_{g_i,0} \} $$
      or equivalently:
      $$ \left(\frac{p}{d}-1, \frac{p}{d} \right)  \cap^{\#} \{ \gamma_d^{*n} \} \geq\sum_{i=1}^r  \left( \frac{p}{d}-1 , \frac{p}{d}\right) \cap^{\#} \{ \sigma_{g_i,0} \}.  $$
\end{thm}

\begin{proof}

    Theorem \ref{th7.1} tells us the following sum in $\mathbb{Z}[(\mathbb{Q}]$ is effective:
   $$ \gamma_d^{*(n+1)} - \left( \sum_{i=1}^r\sum_{j=1}^{\mu_i}\sum_{k=1}^{d-1} \bigg[\frac{\lfloor d\lambda_{ij} \rfloor + k}{d} \bigg] + \sum_{\substack{i=1 \\ \\ d\lambda_{ij}}}^r\sum_{\substack{j=1 \\ \\ \notin \mathbb{Z}}}^{\mu_i} \bigg[\frac{\lfloor d\lambda_{ij} \rfloor  }{d} +1 \bigg] \right) $$
    This is equivalent to a set of statements for every coefficient of $\left[ \frac{p}{d} \right]$ in $\gamma_d^{(n+1)}$ and the coefficient on $\left[\frac{p}{d}\right]$ in the summation to the right of it. That is, for every $p \in \mathbb{Z}$, we have:
    
    \begin{multline*}
        \left\{ \frac{p}{d}  \right\} \cap^{\#} \{\gamma_d^{*(n+1)} \} \geq \sum_{i=1}^r  \sum_{j=1}^{\mu_i} \# \left\{[\lambda{ij}]: \frac{p}{d}=  \frac{\lfloor d\lambda_{ij} \rfloor}{d} + 1  , d\lambda_{ij} \notin \mathbb{Z} \right\}\\
        + \sum_{i=1}^r\sum_{j=1}^{\mu_i}\sum_{k=1}^{d-1} \# \left\{[\lambda_{ij}]: \frac{p}{d} = \frac{\lfloor d\lambda_{ij} \rfloor + k }{d} \right\}
    \end{multline*}

     $$= \sum_{i=1}^r \left\{[\beta]: p-d < d\beta < p \right\} \cap^{\#} \{\sigma_{g_i,0}\} = \sum_{i=1}^r \left( \frac{p}{d}-1 , \frac{p}{d}\right)  \cap^{\#} \{\sigma_{g_i,0} \} $$
    giving our result. The alternative statement follows immediately from the above lemma.
\end{proof}
    
We have now deduced that the bounding argument resulting from our above Theorem \ref{thm1} mimics the form of a particular case of the Varchenko bound. We note that while Varchenko's bound is stronger than our bound, as proven in the following theorem, it is not the case that Varchenko's bounding argument implies the full scope of Theorem \ref{thm1} itself. This is due to the fact that we sacrificed many of the structurally important components of the spectra in order to make this bounding argument in the first place.

\begin{thm}

    Varchenko's bound Theorem \ref{th3.1} implies the formula given in Theorem \ref{th7.7}.

\end{thm}

\begin{proof}

    Let $Z \subseteq \mathbb{P}^n$ be a hypersurface of degree $d$, with only isolated singular points $P_1, \ldots , P_r$. Let $g_i:(\mathbb{C}^n,0) \rightarrow (\mathbb{C},0)$, for $1\leq i \leq r$ denote the corresponding germs defined locally about $P_i$.
    
    Then Theorem \ref{th3.1} tells us that for any $\alpha$, we must have:
    
    $$  (\alpha, \alpha +1)  \cap^{\#} \{ \gamma_{d}^{*n} \} \geq \sum_{i=1}^r  (\alpha, \alpha +1)  \cap^{\#} \{\sigma_{g_i,0} \}  $$
    For any $p \in \mathbb{Z}$, let $\alpha = \frac{p}{d}-1$. Since this can be done for each $p$ the result immediately follows. 
\end{proof}

It is our experience that picking $\alpha \neq \frac{p}{d}$ for some $p \in \mathbb{Z}$, always gives a worse bound than either $p=\lfloor d \alpha \rfloor$ or $p=\lfloor d \alpha \rfloor + 1$. However, upon further inspection of the work of \cite{St85} and \cite{vS20}, we can remedy this gap by adding a set of discrete, and often extraneous, bounds.

The proof of a ``stronger'' bound is actually relatively simple upon a close inspection of the methods detailed in von Straten's exposition \cite{vS20} and a semicontinuity property detailed in \cite{St85}. In particular, the reasoning behind why Varchenko's bound allows for values of $\alpha \in \mathbb{R}$, or really any value $\alpha \neq \frac{p}{d}, p \in \mathbb{Z}$ has more to do with explaining semicontinuity rather than an effort to generate more useful bounds, especially in the case where $f$ is a homogeneous polynomial. 

According to \cite{vS20}, the extension of the semicontinuity to intervals of the form $(\alpha, \alpha + 1 ]$ relies on describing an adjacency of $f + z^m$, which is not necessarily homogeneous, and only extends to new values of $\alpha$ precisely when $m\neq d$. 

We will use the following fact from \cite{St85}

\begin{thm}[Steenbrink 1985]\label{th7} Let $f: \mathbb{C}^{n+1} \rightarrow \mathbb{C}$ have an isolated singularity at $0$, and suppose there exists an adjacency $f \rightsquigarrow g_1 + \ldots + g_r$. Then for every $\alpha \in \mathbb{R}$, $(\alpha, \alpha+1]$ is a semicontinuity domain. In particular, if $f$ is homogeneous of degree $d$, 

 $$  (\alpha, \alpha +1]  \cap^{\#} \{ \gamma_{d}^{*n} \} \geq \sum_{i=1}^r  (\alpha, \alpha +1]  \cap^{\#} \{\sigma_{g_i,0} \}  $$

\end{thm}

We note that \cite{vS20} details how Varchenko's original proof this semicontinuity domain utilizes the Bruce deformation to describe the adjacency between a homogenous polynomial $f$ of degree $d$ and the same local singularities $g_i$ as we describe above. Therefore this inequality, although proven independently, still applies directly to our setting. 

We arrive at the following bound, which we show is stronger than that of Varchenko, and is still discrete in nature. It is merely a combination of our bound and that of Steenbrink when $\alpha = \frac{p}{d}-1$

\pagebreak

\begin{thm}[Open and Half Closed Bound] \label{th7.10}
    Let $f \in \mathbb{C}[x_0, \ldots , x_n]$ be a homogeneous polynomial of degree $d$, and let $\widetilde{V}(f) \subseteq \mathbb{P}^n$ have only $r$ isolated singularities given locally by  $g_i:(\mathbb{C}^n,0)\rightarrow (\mathbb{C},0)$ having corresponding Milnor numbers $\mu_i$ for $1\leq i \leq r$. Let $\sigma_{g_i,0}$ be their corresponding spectra. Then for every $p \in \mathbb{Z}$ we must have:
    
      $$ \left(\frac{p}{d}-1, \frac{p}{d} \right)  \cap^{\#} \{ \gamma_d^{*n} \} \geq\sum_{i=1}^r  \left( \frac{p}{d}-1 , \frac{p}{d}\right) \cap^{\#} \{ \sigma_{g_i,0} \}.  $$

      and 

    $$ \left(\frac{p}{d}-1, \frac{p}{d} \right]  \cap^{\#} \{ \gamma_d^{*n} \} \geq\sum_{i=1}^r  \left( \frac{p}{d}-1 , \frac{p}{d}\right] \cap^{\#} \{ \sigma_{g_i,0} \}.  $$
      
\end{thm}

\begin{thm}

    Our bound, given in Theorem \ref{th7.10} implies the bound given by Varchenko in Theorem \ref{th3.1}.

\end{thm}

\begin{proof}

    Theorem \ref{th3.1} tells us that for any $\alpha$, we must have:
    
    $$  (\alpha, \alpha +1)  \cap^{\#} \{ \gamma_{d}^{*n} \} \geq \sum_{i=1}^r  (\alpha, \alpha +1)  \cap^{\#} \{\sigma_{g_i,0} \}  $$

    The Case $\alpha = \frac{p}{d}-1$  For any $p \in \mathbb{Z}$ is immediately implied by the first bound in our theorem. 

    For $\alpha = \frac{p}{d} -1 + \varepsilon$ for 
    $0<\varepsilon < \frac{1}{d}$, we have:

    $$  \left(\frac{p}{d} -1 + \varepsilon, \frac{p}{d} + \varepsilon \right)  \cap^{\#} \{ \gamma_{d}^{*n} \} =  \left(\frac{p}{d} -1, \frac{p}{d} \right]  \cap^{\#} \{ \gamma_{d}^{*n} \}$$
    and:

    $$\sum_{i=1}^r   \left(\frac{p}{d} -1 + \varepsilon, \frac{p}{d} + \varepsilon \right)  \cap^{\#} \{\sigma_{g_i,0} \} \leq \sum_{i=1}^r   \left(\frac{p}{d} -1 , \frac{p}{d} \right]  \cap^{\#} \{\sigma_{g_i,0} \} $$
    So the latter bound in our theorem implies:

    $$\sum_{i=1}^r   \left(\frac{p}{d} -1 + \varepsilon, \frac{p}{d} + \varepsilon \right)  \cap^{\#} \{\sigma_{g_i,0} \} \leq  \left(\frac{p}{d} -1 + \varepsilon, \frac{p}{d} + \varepsilon \right)  \cap^{\#} \{ \gamma_{d}^{*n} \} $$
    Which covers all of the remaining cases for $\alpha$

\end{proof}

We will now demonstrate the usefulness of Theorem \ref{th7.7} with several examples, which duplicate or improve known results. In the paper \cite{Al03}, the author gives an explicit example of a projective hypersurface $X \subset \mathbb{P}^4$ of degree $3$ with a single $A_{11}$ singularity. In the paper \cite{KL2}, the authors prove that $m=11$ is in fact the maximal $A_m$ singularity that can be present in such a hypersurface with only isolated singularities. The following example illustrates how one can use Theorem \ref{th7.7} to provide another proof that this is in fact the case:

\begin{exmp}

    Let $X \subset \mathbb{P}^4$ be a hypersurface of degree $3$ with only $r$ isolated singularities. Here we have $n=4$ and $d=3$. Assume $X$ has an $A_m$ singularity. Without loss of generality, assume the local equations $g_i$ $1\leq i \leq r$ corresponding to the r singularities are indexed in such a way that $g_1$ corresponds to the $A_m$ singularity. 

    The local normal form of $g_1$ is given by the equation: $z_0^{m+1}+z_1^2+z_2^2+z_3^2$, and the corresponding spectrum is given by:

    $$ \sigma_{g_1,0}= \sum_{j=1}^{\mu_1} [\lambda_{1j}] = \sum_{j=1}^m \left[\frac{j}{m+1} + \frac{1}{2}\right] $$  
    Therefore by Theorem \ref{th7.7}, 
    
    $$1 = \left\{ \frac{2}{3}  \right\} \cap^{\#} \{\gamma_3^{*5} \}  \geq \left( -\frac{2}{6} , \frac{4}{6}\right) \cap^{\#} \{ \sigma_{g_1,0} \} = \# \left\{j : 1\leq j \leq m, \quad j < \frac{m+1}{6}\right\} $$
    However, the last count is greater than 1 whenever $m>11$, and so we must have $m\leq11$. 
    
    We can easily extend this argument to general $n,d$. We get that:
    
        $$ m \leq \frac{2d+(d-2)(n+1)}{2d-(d-2)(n+1)}. $$

\end{exmp}

We now show how Theorem \ref{th7.7} can be used to extend and improve the arguments of Proposition \ref{prop2.2}, and Proposition \ref{prop2.3}. The Hodge-theoretic bounds rely on the assumption that $n+1$ is even to work. Furthermore, one can verify that the former bound is equivalent to using Theorem \ref{th7.7} with $p=d$, which can be improved if we allow the interval some flexibility. 

\begin{prop}

     Let $f \in \mathbb{C}[x_0, x_1, x_2, x_3]$ be homogeneous polynomial of degree $d > 3$, and let $\widetilde{V}(f) \subseteq \mathbb{P}^3$ have only $n_6$ isolated $\widetilde{E_6}$, $n_7$ isolated $\widetilde{E_7}$, and $n_8$ isolated $\widetilde{E_8}$ singularities. Define:
     
     $$ b(d,p) = \binom{3d-p-1}{3} - 4\binom{2d-p}{3} + 6\binom{d-p+1}{3}  $$
     $$ =\frac{d^3}{6} + \frac{d^2p}{2} - \frac{p^3}{2} + \frac{dp^2}{2} + -d^2 -2dp +p^2   +\frac{11d}{6} + \frac{p}{2}  - 1$$
     
     Then:
     
     \begin{enumerate}
         \item $7n_6 \leq b(d,p) \leq \frac{31}{54}d^3 - \frac{13}{18}d^2 +4d + \frac{1}{2}$ \hspace{3mm} for \hspace{3mm} $p = \left\lfloor \frac{2d}{3} \right\rfloor +1$
         \item $7n_6 + 8n_7 \leq b(d,p) \leq \frac{235}{384}d^3 - \frac{11}{16}d^2 +\frac{101}{24}d + \frac{1}{2} $ \hspace{3mm} for \hspace{3mm} $p = \left\lfloor \frac{3d}{4} \right\rfloor +1$
         \item $7n_6 + 8n_7 +9n_8 \leq b(d,p) \leq \frac{277}{432}d^3 - \frac{23}{36}d^2 +\frac{53}{12}d + \frac{1}{2} $ \hspace{3mm} for \hspace{3mm} $p = \left\lfloor \frac{5d}{6} \right\rfloor +1$
         
     \end{enumerate}

\end{prop}

\begin{proof}
   
   We note that the normal forms and spectra for each singularity type are as follows:
   
   \begin{enumerate}
       \item $\widetilde{E_6}$:  $\quad x^3+y^3+z^3$ $ \quad \sigma_{\widetilde{E_6}} = [0] +3 \left[ \frac{1}{3}\right] + 3 \left[ \frac{2}{3}\right] + [1]$
       \item $\widetilde{E_7}$:  $\quad x^2+y^4+z^4$ $ \quad \sigma_{\widetilde{E_7}} = [0] +2 \left[ \frac{1}{4}\right] + 3 \left[ \frac{1}{2}\right] + 2\left[ \frac{3}{4}\right] + [1]$
       \item $\widetilde{E_8}$:  $\quad x^2+y^3+z^6$ $ \quad \sigma_{\widetilde{E_8}} = [0] +1 \left[ \frac{1}{6}\right] + 2 \left[ \frac{1}{3}\right] + 2\left[ \frac{1}{2}\right] + 2\left[ \frac{2}{3}\right] + 1\left[ \frac{5}{6}\right] + [1]$
   \end{enumerate}
   
   We apply Theorem \ref{th7.7} for the choices of $p= \left\lfloor \frac{2d}{3} \right\rfloor +1, \left\lfloor \frac{3d}{4} \right\rfloor +1$, and $\left\lfloor \frac{5d}{6} \right\rfloor +1$, respectively. The $b(d, p)$ are simply calculated using Proposition \ref{prop7.4}.  Since these values of $p$ are increasing, and the values of $\frac{p}{d}- 1 < 0 $ in any of these choices, our bound for $8n_7$, duplicates as a bound for $7n_6 + 8n_7$ and our bound for $9n_8$, duplicates as a bound for $7n_6 + 8n_7 + 9n_8$. The polynomials bounding $b(d,p)$, are obtained from inequalities of the form $ \frac{2d}{3} \leq \left\lfloor \frac{2d}{3} \right\rfloor +1 \leq \frac{2d}{3} + 1 $, and plugging these values into the polynomial form of $b(d,p)$, depending on the sign of $p^m$ in each summand.
   
\end{proof}   
   
\begin{exmp}
    In particular, we can compare the first bound with that of Proposition 2.2. The bounds for the following manually calculated bounds and those of Proposition 2.2 for $4\leq d \leq 9$
    
    \begin{center}
        \begin{tabular}{||c | c | c | c | c | c | c ||}
        \hline 
        &&&&&& \\
        bound/d & 4 & 5 & 6 & 7 & 8 & 9 \\
        \hline
        &&&&&& \\
        $ \frac{1}{7} b(d,p) $ & 2 & 5 & 11 & 17 & 29 & 45 \\
        \hline
        &&&&&& \\
        $ \frac{1}{6}( h^{1,1}_{3,d} -1 ) $ & 3 & 7 & 14 & 24 & 38 & 56 \\
        \hline

        \end{tabular}
    \end{center}
    for $d\geq 10$, the values $ \frac{1}{6}( h^{1,1}_{3,d} -1 ) \geq  \frac{1}{7}\left(\frac{31}{54}d^3 - \frac{13}{18}d^2 +4d + \frac{1}{2} \right) \geq \frac{1}{7}b(d,p) $ for $p = \left\lfloor \frac{2d}{3} \right\rfloor +1$
    So this bound is always better than the hodge theoretic bound given in Proposition 2.2
\end{exmp}    

\begin{prop}

    Let $f \in \mathbb{C}[x_0, x_1, x_2. x_3, x_4]$ be homogeneous polynomial of degree $d > 2$, and let $\widetilde{V}(f) \subseteq \mathbb{P}^4$ have only isolated $A_{2m+1}$ singularities. Then the number $r$ of singular points is bounded by
    
    $$ r \leq  \begin{cases} 
     \frac{1}{2m+1}\left[ \frac{115}{192}d^4 -\frac{115}{48}d^3+ \frac{185}{48}d^2 - \frac{35}{12}d + 1 \right]  & d\equiv 0 \mod 2 \\
      \frac{1}{2m+1}\left[ \frac{115}{192}d^4 -\frac{115}{48}d^3 + \frac{355}{96}d^2 - \frac{125}{48}d + \frac{45}{64} \right]  & d\equiv 1 \mod 2, d>m+1\\
   \end{cases}
    $$

\end{prop}

\begin{proof}

    Let $X \subset \mathbb{P}^4$ be a hypersurface of degree $d$ with only $r$ isolated $A_{2m+1}$ singularities Let the local equations be given by $g_i$ $1\leq i \leq r$. 

    The local normal form of $g_i$ for $1\leq r \leq$ is given by the equation: $z_0^{2m+2}+z_1^2+z_2^2+z_3^2$, and the corresponding spectra are given by:

    $$ \sigma_{g_i,0}= \sum_{j=1}^{\mu_1} [\lambda_{ij}] = \sum_{j=1}^{2m+1} \left[\frac{j}{2m+2} + \frac{1}{2}\right] $$
    Pick $p= \frac{3d}{2} $ for $d$ even, and $p=\frac{3d+1}{2}$ for $d$ odd (assuming $d>m+1$).  Then by \ref{th7.7}, 
    
    $$\left\{ \frac{p}{d}  \right\} \cap^{\#} \{\gamma_d^{*5} \}  \geq \sum_{i=1}^r \left( \frac{p}{d} - 1 , \frac{p}{d}\right) \cap^{\#} \{ \sigma_{g_i,0} \} = r(2m+1) $$
    
    A simple calculation of the left hand side using Proposition \ref{prop7.4} for $d\equiv 0,1 \mod2$ gives the desired result.
\end{proof}

We note, however, that this bound is best for larger $m$. One can verify, for example, that our bound on the number of $A_1$ singularities given above does better (this is the case where $m=0$). This is because the bounds above were chosen to give a convenient polynomial bound that does best for all $m$. Technically, when working with a particular $m$, a better bound can be found by choosing $p=\left\lfloor \frac{(3m+2)d}{2m+2} \right\rfloor +1 $. In particular this can be done even for the cases when $d$ is odd and $m+1 \geq d$.

\section{Appendix: comparing the bounds}\label{S8}

There is in fact quite a long history of bounding the number of singularities of projective hypersurfaces, most notably bounding the number of nodes (also called an $A_1$ singularity or ordinary double point).  In the tables that follow, ``naive'' denotes the vanishing cycle sequence method of Section \ref{S2}.  Note that it makes no prediction for nodes on a threefold (second table).

\pagebreak

$A_1$ singularities, $n=3$:
\begin{center}
 \begin{tabular}{||c | c | c | c | c||}
     \hline
      d & Naive & Eigenvalue & Conical & Sharp\\ [.5ex] 
     \hline \hline
     &&&&\\
     Eq & $\frac{2}{3}d^3-2d^2+\frac{7}{3}d-1$ & $d^3-4d^2+6d-4$ & $\frac{23}{48}d^3-\frac{9}{8}d^2 +\frac{5}{6}d$, even d & \\ 
    &&&& \\
        &  &  & $\frac{23}{48}d^3-\frac{23}{16}d^2 +\frac{78}{48}d +\frac{9}{16}$ odd  & \\[1ex]
     \hline
     1 & 0 &  &  & 0 \\ 
     \hline
     2 & 1 &  &  & 1 \\
     \hline
     3 & 6 & 5 & 4 & 4 \\
     \hline
     4 & 19 & 20 & 16 & 16 \\
     \hline
     5 & 44 & 51 & 31 & 31\\
     \hline
     6 & 85 & 104 & 68 & 65\\
    \hline
     7 & 146 & 185 & 104 & 99-104\\
    \hline
     10 & 489 & 656 & 375  & \\
    \hline
     20 & 4579 & 6516 & 3400 & \\
    \hline
     30 & 16269 & 23576 & 11950 & \\
    \hline
     40 & 39559 & 57836 & 28900 & \\
     \hline
     50 & 78449 & 115296 & 57125 & \\
     \hline
     100 & 646899 & 960596 & 468000 & \\
    \hline
     1,000 & 664668999 & 996005996 & 478042500 & \\ [1ex]
     \hline
\end{tabular}
\end{center}

$A_1$ singularities, $n=4$:
\begin{center}
 \begin{tabular}{||c | c | c | c | c||} 
     \hline
      d & Eigenvalue & Conical & Sharp\\ [0.5ex] 
     \hline \hline
     &&& \\
     Eq   & $d^4-5d^3+10d^2-10d+5$  & $\frac{11}{24}d^4-\frac{19}{12}d^3+\frac{49}{24}d^2-\frac{11}{12}d$ & \\ [1ex]
     \hline
     1 &  &  & 0  \\ 
     \hline
     2 &  &  & 1  \\
     \hline
     3 & 11 & 10 & 10  \\
     \hline
     4  & 61 & 45 & 45\\
     \hline
     5   & 205 & 135 & 130-135 \\
     \hline
     6   & 521 & 320 & \\
    \hline
     7   & 1111 & 651 & \\
    \hline
     10   & 5905 & 3195 & \\
    \hline
     20   & 123805 & 61465 & \\
    \hline
     30   & 683705 & 330310 & \\
    \hline
     40   & 2255605 & 1075230 & \\
     \hline
     50   & 5649505 & 2671725  & \\
     \hline
     100   & 95099005 & 44270325 & \\ [1ex]
     \hline
\end{tabular}
\end{center}

\pagebreak

$\widetilde{E_6}$ singularities, $n=3$:
\begin{center}
 \begin{tabular}{||c | c | c | c | c||}
     \hline
      d & Naive & Eigenvalue & Conical & Sharp\\ [.5ex] 
     \hline \hline
     &&&&\\
     Eq & $\frac{1}{9}d^3-\frac{1}{3}d^2+\frac{7}{18}d-\frac{1}{6}$ & $\frac{1}{8}d^3-\frac{1}{2}d^2+\frac{3}{4}d-\frac{1}{2}$ & $\frac{1}{7}b(d,p)$ & \\ [1ex]
     \hline
     1 & 0 &  &  & 0 \\ 
     \hline
     2 & 0 &  &  & 0 \\
     \hline
     3 & 1 &  &  & 1 \\
     \hline
     4 & 3 & 3 & 2 & 1 \\
     \hline
     5 & 7 & 6 & 5 & \\
     \hline
     6 & 14 & 13 & 11 & \\
    \hline
     7 & 24 & 23 & 17 & \\
    \hline
     10 & 82 & 82 & 60  & \\
    \hline
     20 & 763 & 815 & 570 & \\
    \hline
     30 & 2712 & 2947 & 2040 & \\
    \hline
     40 & 6593 & 7230 & 4865 & \\
     \hline
     50 & 13075 & 14412 & 9706 & \\
     \hline
     100 & 107817 & 120075 & 79577 & \\
    \hline
     1,000 & 110778167 & 124500750 & 81764819 & \\ [1ex]
     \hline
\end{tabular}
\end{center}

\pagebreak

For reference, we give the output of $\gamma_d^{*(n+1)}$ for select values:

\begin{center}
    \begin{tabular}{||c | c | c ||} 
         \hline
         & & \\
         n & d & $\gamma_d^{*(n+1)}$ \\ [0.5ex] 
         \hline\hline
         &&\\
         2 & 2 & $1\left[\frac{1}{2}\right]$ \\[1ex]
         \hline
         &&\\
         2 & 3 & $1 [0] + 3\left[\frac{1}{3}\right]+ 3\left[\frac{2}{3}\right] +1[1]$ \\[1ex]
         \hline
         &&\\
         2 & 4 & $1\left[-\frac{1}{4}\right] + 3\left[0\right]+ 6\left[\frac{1}{4}\right] + 7\left[\frac{1}{2}\right] + 6\left[ \frac{3}{4} \right] + 3\left[ 1 \right] + 1\left[\frac{5}{4}\right]$  \\[1ex]
         \hline
         &&\\
         2 & 5 & $1\left[-\frac{2}{5}\right] + 3\left[-\frac{1}{5}\right]+ 6\left[0\right]+10\left[\frac{1}{5}\right]+12\left[\frac{2}{5}\right]$ \\[1ex] & & $+12\left[\frac{3}{5}\right]+10\left[\frac{4}{5}\right]+6\left[1\right]+3\left[\frac{6}{5}\right]+1\left[\frac{7}{5}\right]$ \\[1ex]
         \hline
         &&\\
          2 & 6 & $1\left[-\frac{1}{2}\right] + 3\left[-\frac{1}{3}\right]+ 6\left[-\frac{1}{6}\right]+10\left[0\right]+15\left[\frac{1}{6}\right]+18\left[\frac{1}{3}\right]+19\left[\frac{1}{2}\right]$ \\[1ex] & & $+18\left[\frac{2}{3}\right]+15\left[\frac{5}{6}\right]+10\left[1\right]+6\left[\frac{7}{6}\right]+3\left[\frac{4}{3}\right]+1\left[\frac{3}{2}\right]$ \\[1ex]
         \hline
         &&\\
          2 & 7 & $1\left[-\frac{4}{7}\right] + 3\left[-\frac{3}{7}\right]+ 6\left[-\frac{2}{7}\right]+10\left[-\frac{1}{7}\right]+15\left[0\right]+21\left[\frac{1}{7}\right]+25\left[\frac{2}{7}\right]+27\left[\frac{3}{7}\right]$ \\[1ex] & & $+27\left[\frac{4}{7}\right]+25\left[\frac{5}{7}\right]+21\left[\frac{6}{7}\right]+15\left[1\right]+10\left[\frac{8}{7}\right]+6\left[\frac{9}{7}\right]+3\left[\frac{10}{7}\right]+1\left[\frac{11}{7}\right]$  \\[1ex]
         \hline
         &&\\
          3 & 2 & $1[1]$ \\[1ex]
         \hline
         &&\\
          3 & 3 & $1\left[\frac{1}{3}\right] + 4\left[\frac{2}{3}\right] + 6\left[1\right]+ 4\left[\frac{4}{3}\right]+ 1\left[\frac{5}{3}\right]$ \\[1ex]
         \hline
         &&\\
          3 & 4 & $1[0]+4\left[\frac{1}{4}\right] + 10\left[\frac{1}{2}\right]+ 16\left[\frac{3}{4}\right] + 19\left[1\right] + 16\left[ \frac{5}{4} \right] + 10\left[\frac{3}{2}\right]  +4\left[\frac{7}{4}\right] +  1\left[2\right]$  \\[1ex]
         \hline
         &&\\
          3 & 5 & $1\left[-\frac{1}{5}\right] + 4\left[0\right] +10\left[\frac{1}{5}\right]+ 20\left[\frac{2}{5}\right]+31\left[\frac{3}{5}\right]+40\left[\frac{4}{5}\right]+44\left[1\right]$ \\[1ex] & & $+40\left[\frac{6}{5}\right]+31\left[\frac{7}{5}\right]+20\left[\frac{8}{5}\right]+10\left[\frac{9}{5}\right]+4\left[2\right]+1\left[\frac{11}{5}\right]$ \\[1ex] 
         \hline
         &&\\
          3 & 6 & $1\left[-\frac{1}{3}\right]+ 4\left[-\frac{1}{6}\right]+10\left[0\right]+20\left[\frac{1}{6}\right]+35\left[\frac{1}{3}\right]+52\left[\frac{1}{2}\right]+68\left[\frac{2}{3}\right]+80\left[\frac{5}{6}\right]+85\left[1\right]$ \\[1ex] & & $+80\left[\frac{7}{6}\right]+68\left[\frac{4}{3}\right]+52\left[\frac{3}{2}\right]+35\left[\frac{5}{3}\right]+20\left[\frac{11}{6}\right]+10\left[2\right]+4\left[\frac{13}{6}\right]+1\left[\frac{7}{3}\right]$ \\[1ex]
         \hline
         &&\\
          3 & 7 & $1\left[-\frac{3}{7}\right]+ 4\left[-\frac{2}{7}\right]+10\left[-\frac{1}{7}\right]+20\left[0\right]+35\left[\frac{1}{7}\right]+56\left[\frac{2}{7}\right]+80\left[\frac{3}{7}\right]+104\left[\frac{4}{7}\right]$\\[1ex] & &
          $+125\left[\frac{5}{7}\right]+140\left[\frac{6}{7}\right] +146\left[1\right] +140\left[\frac{8}{7}\right]+125\left[\frac{9}{7}\right]+104\left[\frac{10}{7}\right]+80\left[\frac{11}{7}\right]+56\left[\frac{12}{7}\right]$
          \\[1ex] & &
          $+35\left[\frac{13}{7}\right]+20\left[2\right]+10\left[\frac{15}{7}\right]+4\left[\frac{16}{7}\right] +1\left[\frac{17}{7}\right]$\\[1ex]
         \hline
         \end{tabular}
         
    \begin{tabular}{||c | c | c ||}
        \hline
         & & \\
         n & d & $\gamma_d^{*(n+1)}$ \\ [0.5ex] 
         \hline\hline
          & & \\
          4 & 2 & $\left[\frac{3}{2}\right]$ \\[1ex]
         \hline
         &&\\
          4 & 3 & $1\left[\frac{2}{3}\right]+ 5\left[1\right] + 10\left[\frac{4}{3}\right]+ 10\left[\frac{5}{3}\right]+ 5\left[2\right]+ 1\left[\frac{7}{3}\right]$ \\[1ex]
         \hline
         &&\\
          4 & 4 & $1\left[\frac{1}{4}\right] + 5\left[\frac{1}{2}\right]+ 15\left[\frac{3}{4}\right] + 30\left[1\right] + 45\left[ \frac{5}{4} \right] + 51\left[\frac{3}{2}\right]$ \\[1ex] & &
          $+45\left[\frac{7}{4}\right] +  30\left[2\right]+15\left[\frac{9}{4}\right]+5\left[\frac{5}{2}\right]+1\left[\frac{11}{4}\right]$ \\[1ex]
         \hline
         &&\\
          4 & 5 &  $1\left[0\right] +5\left[\frac{1}{5}\right]+ 15\left[\frac{2}{5}\right]+35\left[\frac{3}{5}\right]+65\left[\frac{4}{5}\right]+101\left[1\right]+135\left[\frac{6}{5}\right]+155\left[\frac{7}{5}\right]$ \\[1ex] & & $+155\left[\frac{8}{5}\right]+135\left[\frac{9}{5}\right]+101\left[2\right]+65\left[\frac{11}{5}\right]+35\left[\frac{12}{5}\right]+15\left[\frac{13}{5}\right]+5\left[\frac{14}{5}\right]+1\left[3\right]$  \\[1ex]
         \hline
         &&\\
          4 & 6 &  $1\left[-\frac{1}{6}\right]+5\left[0\right]+15\left[\frac{1}{6}\right]+35\left[\frac{1}{3}\right]+70\left[\frac{1}{2}\right]+121\left[\frac{2}{3}\right]+185\left[\frac{5}{6}\right]+255\left[1\right]$ \\[1ex] & & $+320\left[\frac{7}{6}\right]+365\left[\frac{4}{3}\right]+381\left[\frac{3}{2}\right]+365\left[\frac{5}{3}\right]+320\left[\frac{11}{6}\right]+255\left[2\right]+185\left[\frac{13}{6}\right]$ \\[1ex] & & $+121\left[\frac{7}{3}\right]+70\left[\frac{5}{2}\right]+35\left[\frac{8}{3}\right]+15\left[\frac{17}{6}\right]+5\left[3\right]+1\left[\frac{19}{6}\right]$  \\[1ex]
         \hline
         &&\\
          4 & 7 &  $1\left[-\frac{2}{7}\right]+5\left[-\frac{1}{7}\right]+15\left[0\right]+35\left[\frac{1}{7}\right]+70\left[\frac{2}{7}\right]+126\left[\frac{3}{7}\right]+205\left[\frac{4}{7}\right]$\\[1ex] & &
          $+305\left[\frac{5}{7}\right]+420\left[\frac{6}{7}\right] +540\left[1\right] +651\left[\frac{8}{7}\right]+735\left[\frac{9}{7}\right]+780\left[\frac{10}{7}\right]+780\left[\frac{11}{7}\right]$\\[1ex] & &
          $+735\left[\frac{12}{7}\right]
          +651\left[\frac{13}{7}\right]+540\left[2\right]+420\left[\frac{15}{7}\right]+305\left[\frac{16}{7}\right] +205\left[\frac{17}{7}\right]+126\left[\frac{18}{7}\right]$\\[1ex] & &
          $+70\left[\frac{19}{7}\right]+35\left[\frac{20}{7}\right]+15\left[3\right]+5\left[\frac{22}{7}\right]+1\left[\frac{23}{7}\right]$ \\[1ex]
         \hline
         &&\\
         5 & 2 & $1\left[2\right]$ \\[1ex]
         \hline
         &&\\
         5 & 3 & $1\left[1\right]+ 6\left[\frac{4}{3}\right]+ 15\left[\frac{5}{3}\right]+ 20\left[2\right]+ 15\left[\frac{7}{3}\right]+ 6\left[\frac{8}{3}\right]+ 1\left[3\right]$  \\[1ex]
         \hline
         &&\\
         5 & 4 & $1\left[\frac{1}{2}\right]+ 6\left[\frac{3}{4}\right] + 21\left[1\right] + 50\left[ \frac{5}{4} \right] + 90\left[\frac{3}{2}\right]+126\left[\frac{7}{4}\right] +  141\left[2\right]$\\[1ex] & &
          $+126\left[\frac{9}{4}\right]+90\left[\frac{5}{2}\right]+50\left[\frac{11}{4}\right]+21\left[3\right]+6\left[\frac{13}{4}\right]+1\left[\frac{7}{2}\right]$ \\[1ex]
         \hline
         &&\\
         5 & 5 &  $1\left[\frac{1}{5}\right]+ 6\left[\frac{2}{5}\right]+21\left[\frac{3}{5}\right]+56\left[\frac{4}{5}\right]+120\left[1\right]+216\left[\frac{6}{5}\right]+336\left[\frac{7}{5}\right]$\\[1ex] & &
          $ +456\left[\frac{8}{5}\right]+546\left[\frac{9}{5}\right]+580\left[2\right]+546\left[\frac{11}{5}\right]+456\left[\frac{12}{5}\right]+336\left[\frac{13}{5}\right]$\\[1ex] & &
          $+216\left[\frac{14}{5}\right]+120\left[3\right]+56\left[\frac{16}{5}\right]+21\left[\frac{17}{5}\right]+6\left[\frac{18}{5}\right]+1\left[\frac{19}{5}\right]$  \\[1ex]
         \hline
         &&\\
         5 & 6 &  $1\left[0\right]+6\left[\frac{1}{6}\right]+21\left[\frac{1}{3}\right]+56\left[\frac{1}{2}\right]+126\left[\frac{2}{3}\right]+246\left[\frac{5}{6}\right]+426\left[1\right]$\\[1ex] & &
          $+666\left[\frac{7}{6}\right]+951\left[\frac{4}{3}\right]+1246\left[\frac{3}{2}\right]+1506\left[\frac{5}{3}\right]+1686\left[\frac{11}{6}\right]+1751\left[2\right]$\\[1ex] & &
          $+1686\left[\frac{13}{6}\right]+1506\left[\frac{7}{3}\right]+1246\left[\frac{5}{2}\right]+951\left[\frac{8}{3}\right]+666\left[\frac{17}{6}\right]+426\left[3\right]$\\[1ex] & &
          $+246\left[\frac{19}{6}\right]+126\left[\frac{10}{3}\right]+56\left[\frac{7}{2}\right]+21\left[\frac{11}{3}\right]+6\left[\frac{23}{6}\right]+1\left[4\right]$  \\[1ex]
         \hline
     \end{tabular}
         
     \begin{tabular}{||c | c | c ||} 
         \hline
         & & \\
         n & d & $\gamma_d^{*(n+1)}$ \\ [0.5ex] 
         \hline\hline
         &&\\
         5 & 7 & $1\left[-\frac{1}{7}\right]+6\left[0\right]+21\left[\frac{1}{7}\right]+56\left[\frac{2}{7}\right]+126\left[\frac{3}{7}\right]+252\left[\frac{4}{7}\right]+456\left[\frac{5}{7}\right]$\\[1ex] & &
          $+756\left[\frac{6}{7}\right] +1161\left[1\right] +1666\left[\frac{8}{7}\right]+2247\left[\frac{9}{7}\right]+2856\left[\frac{10}{7}\right]+3431\left[\frac{11}{7}\right]$\\[1ex] & &
          $+3906\left[\frac{12}{7}\right]
          +4221\left[\frac{13}{7}\right]+4332\left[2\right]+4221\left[\frac{15}{7}\right]+3906\left[\frac{16}{7}\right] +3431\left[\frac{17}{7}\right]$\\[1ex] & &
          $+2856\left[\frac{18}{7}\right]
          +2247\left[\frac{19}{7}\right]+1666\left[\frac{20}{7}\right]+1161\left[3\right]+756\left[\frac{22}{7}\right]+456\left[\frac{23}{7}\right]$\\[1ex] & &
          $+252\left[\frac{24}{7}\right]+126\left[\frac{25}{7}\right]+56\left[\frac{26}{7}\right]+21\left[\frac{27}{7}\right]+6\left[4\right]+1\left[\frac{29}{7}\right]$ \\[1ex]
         \hline
         &&\\
         6 & 2 & $1\left[\frac{5}{2}\right]$  \\[1ex]
         \hline
         &&\\
         6 & 3 & $1\left[\frac{4}{3}\right]+ 7\left[\frac{5}{3}\right]+ 21\left[2\right]+ 35\left[\frac{7}{3}\right]+ 35\left[\frac{8}{3}\right]+ 21\left[3\right]+ 7\left[\frac{10}{3}\right]+ 1\left[\frac{11}{3}\right]$  \\[1ex]
         \hline
         &&\\
         6 & 4 & $1\left[\frac{3}{4}\right] + 7\left[1\right] + 28\left[ \frac{5}{4} \right] + 77\left[\frac{3}{2}\right]+161\left[\frac{7}{4}\right] +  266\left[2\right]+357\left[\frac{9}{4}\right]+393\left[\frac{5}{2}\right]$\\[1ex] & &
          $+357\left[\frac{11}{4}\right]+266\left[3\right]+161\left[\frac{13}{4}\right]+77\left[\frac{7}{2}\right]+28\left[\frac{15}{4}\right]+7\left[\frac{8}{2}\right]+1\left[\frac{17}{4}\right]$  \\[1ex]
         \hline
         &&\\
         6 & 5 & $1\left[\frac{2}{5}\right]+7\left[\frac{3}{5}\right]+28\left[\frac{4}{5}\right]+84\left[1\right]+203\left[\frac{6}{5}\right]+413\left[\frac{7}{5}\right] +728\left[\frac{8}{5}\right]+1128\left[\frac{9}{5}\right]$\\[1ex] & &
          $+1554\left[2\right]+1918\left[\frac{11}{5}\right]+2128\left[\frac{12}{5}\right]+2128\left[\frac{13}{5}\right]+1918\left[\frac{14}{5}\right]+1554\left[3\right]$\\[1ex] & &
          $+1128\left[\frac{16}{5}\right]+728\left[\frac{17}{5}\right]+413\left[\frac{18}{5}\right]+203\left[\frac{19}{5}\right]+84\left[4\right]+28\left[\frac{21}{5}\right]+7\left[\frac{22}{5}\right]+1\left[\frac{23}{5}\right]$  \\[1ex]
         \hline
         &&\\
         6 & 6 & $1\left[\frac{1}{6}\right]+7\left[\frac{1}{3}\right]+28\left[\frac{1}{2}\right]+84\left[\frac{2}{3}\right]+210\left[\frac{5}{6}\right]+455\left[1\right]+875\left[\frac{7}{6}\right]$\\[1ex] & &
          $+1520\left[\frac{4}{3}\right]+2415\left[\frac{3}{2}\right]+3535\left[\frac{5}{3}\right]+4795\left[\frac{11}{6}\right]+6055\left[2\right]+7140\left[\frac{13}{6}\right]+7875\left[\frac{7}{3}\right]$\\[1ex] & &
          $+8135\left[\frac{5}{2}\right]+7875\left[\frac{8}{3}\right]+7140\left[\frac{17}{6}\right]+6055\left[3\right]
          +4795\left[\frac{19}{6}\right]+3535\left[\frac{10}{3}\right]$\\[1ex] & &
          $+2415\left[\frac{7}{2}\right]+1520\left[\frac{11}{3}\right]+875\left[\frac{23}{6}\right]+210\left[4\right]+84\left[\frac{25}{6}\right]+28\left[\frac{13}{3}\right]+7\left[\frac{9}{2}\right]+1\left[\frac{29}{6}\right]$   \\[1ex]
         \hline
         &&\\
         6 & 7 & $1\left[0\right]+7\left[\frac{1}{7}\right]+28\left[\frac{2}{7}\right]+84\left[\frac{3}{7}\right]+210\left[\frac{4}{7}\right]+462\left[\frac{5}{7}\right]+917\left[\frac{6}{7}\right]$\\[1ex] & &
          $ +1667\left[1\right] +2807\left[\frac{8}{7}\right]+4417\left[\frac{9}{7}\right]+6538\left[\frac{10}{7}\right]+9142\left[\frac{11}{7}\right]+12117\left[\frac{12}{7}\right]$\\[1ex] & &
          $
          +15267\left[\frac{13}{7}\right]+18327\left[2\right]+20993\left[\frac{15}{7}\right]+22967\left[\frac{16}{7}\right] +24017\left[\frac{17}{7}\right]+24017\left[\frac{18}{7}\right]$\\[1ex] & &
          $
          +22967\left[\frac{19}{7}\right]+20993\left[\frac{20}{7}\right]+18327\left[3\right]+15267\left[\frac{22}{7}\right]+12117\left[\frac{23}{7}\right]+9142\left[\frac{24}{7}\right]$\\[1ex] & &
          $+6538\left[\frac{25}{7}\right]+4417\left[\frac{26}{7}\right]+2807\left[\frac{27}{7}\right]+1667\left[4\right]+917\left[\frac{29}{7}\right]+462\left[\frac{30}{7}\right]+210\left[\frac{31}{7}\right]$\\[1ex] & &
          $+84\left[\frac{32}{7}\right]+28\left[\frac{33}{7}\right]+7\left[\frac{34}{7}\right]+1\left[5\right]$   \\[1ex]
         \hline
    \end{tabular}
\end{center}

\pagebreak

\end{document}